\documentclass[10pt]{article}
\usepackage[ansinew]{inputenc}
\usepackage[french]{babel}
\usepackage{amsmath, amsthm, amssymb}
\usepackage{hyperref}
\hypersetup{pdfpagemode=FullScreen,colorlinks=true}

\title{Pincement du plan hyperbolique complexe}
\author{Pierre Pansu$^{1,2}$\footnote{$^{1}$ Univ Paris-Sud, Laboratoire de Mathématiques d'Orsay, Orsay, F-91405 ;\hfill\eject\indent\hskip 6.3pt $^{2}$ CNRS, Orsay, F-91405.}
}

\date{\today}

\newtheorem{theo}{Th{\'e}or{\`e}me}
\newtheorem{lemme}{Lemme}
\newtheorem*{stheo}{Th{\'e}or{\`e}me}
\newtheorem{prop}[lemme]{Proposition}
\newtheorem{cor}[lemme]{Corollaire}
\newtheorem{defi}[lemme]{D{\'e}finition}
\newtheorem{rem}[lemme]{Remarque}

\def\preuve{\par\medskip\noindent {\bf Preuve.}{\hskip1em}}
\def\remarque{\par\medskip\noindent {\bf Remarque.}{\hskip1em}}
\def\exemple{\par\medskip\noindent {\bf Exemple.}{\hskip1em}}
\def\qed{~q.e.d.}

\def\R{{\bf R}}
\def\C{{\bf C}}
\def\Z{{\bf Z}}
\def\n#1{{\parallel #1 \parallel}}
\def\l{\lambda}
\def\f{\mathfrak{d}}

\def\L{\Lambda}
\def\N{\mathfrak{n}}
\def\H{\mathfrak{h}}
\def\G{\mathfrak{g}}

\def\dim{{\rm dim}\,}
\def\ker{{\rm ker}\,}
\def\im{{\rm im}\,}
\def\tr{{\rm tr}\,}

\def\ddt{\frac{\partial}{\partial t}}
\def\con{\hbox{const.}}
\def\d{\displaystyle}
\def\inv#1{\frac{1}{#1}}
\def\op#1{\Omega^{#1,p}}

\def\psip#1{\Psi^{#1,p}}

\def\lpo#1{L^{p}\Omega^{#1}}

\def\lmpo#1{L_{-1}^{p}\Omega^{#1}}

\def\omp#1{\Omega_{-1}^{#1,p}}
\def\hp#1{H^{#1,p}}

\def\tp#1{T^{#1,p}}

\def\rp#1{R^{#1,p}}

\def\bp#1{{\cal B}^{#1,p}}
\def\p#1#2{{\bf p}_{#1}^{#2}}
\def\lp{_{L^{p}}}

\def\ali{\medskip}
\def\opp#1{\Gamma(\L^{#1}_{+(p)})}
\def\opm#1{\Gamma(\L^{#1}_{-(p)})}
\def\dplus{d_{+}}
\def\dmoins{d_{-}}

\begin{document}
\maketitle
\begin{quote}
{\small
RESUME. On montre que la cohomologie $L^p$ des espaces symétriques de rang un, de type non compact, est séparée dans des intervalles ou cela ne résulte pas seulement du pincement de la courbure. En utilisant la structure multiplicative sur la cohomologie $L^p$, on montre que le plan hyperbolique complexe n'est pas quasiisométrique à une variété riemannienne simplement connexe à courbure sectionnelle strictement $-\frac{1}{4}$-pincée. Malheureusement, la méthode ne s'étend pas aux autres espaces symétriques de rang un.

ABSTRACT. $L^p$-cohomology of rank one symmetric spaces of noncompact type is shown to be Hausdorff for values of $p$ where this does not follow from curvature pinching. Using the multiplicative structure on $L^p$-cohomology, it is shown that no simply connected Riemannian manifold with strictly $-\frac{1}{4}$-pinched sectional curvature can be quasiisometric to complex hyperbolic plane. Unfortunately, the method does not extend to other rank one symmetric spaces.
}
\end{quote}

\section{Introduction}
\label{intro}

\subsection{Pincement négatif}

Si $-1\leq\delta<0$, on dit qu'une variété riemannienne est {\em $\delta$-pinc\'ee} s'il existe $a>0$ tel que sa courbure sectionnelle soit comprise entre $-a$ et $\delta a$. 

Par exemple, l'espace hyperbolique réel est $-1$-pincé. Les espaces
sym\'etriques de rang un de type non compact \`a courbure non
constante sont $-\frac{1}{4}$-pinc{\'e}s. Il s'agit des espaces hyperboliques complexes $\mathbb{C}H^{m}$, $m\geq 2$, des espaces hyperboliques quaternioniens $\mathbb{H}H^{m}$, $m\geq 2$, et du plan hyperbolique des octaves de Cayley $\mathbb{O}H^{2}$.

\bigskip

Sur ces variétés, peut-on changer de métrique riemannienne et améliorer le pincement ? On demande à la nouvelle métrique de rester équivalente à la précédente. Cette question élégante, posée par M. Gromov, \cite{Gromovasympt}, a fait l'objet d'une tentative infructueuse, \cite{Ppince}, utilisant la torsion en cohomologie $L^p$. Dans le présent article, on complète \cite{Ppince}, en montrant que la torsion en cohomologie $L^p$ des espaces symétriques de rang un s'annule là où l'on aurait préféré qu'elle ne s'annule pas. En outre, en utilisant le cup-produit, on résoud le problème dans le cas du plan hyperbolique complexe, mais, malheureusement, seulement dans ce cas là.

\subsection{Cohomologie \texorpdfstring{$L^{p}$}{}}
\label{defcohomologie}
Soit $M$ une vari{\'e}t{\'e} riemannienne. Soit $p>1$ un r{\'e}el. On note
$\lpo{*}(M)$ l'espace de Banach des formes diff{\'e}rentielles $L^p$ et $\Omega^{*,p}(M)=\lpo{*}(M)\cap d^{-1}\lpo{*}(M)$ l'espace des formes diff{\'e}rentielles $L^p$ dont la diff{\'e}rentielle ext{\'e}rieure est aussi $L^p$, muni de la norme
\begin{eqnarray*}
\n{\omega}_{\op{*}}=\left(\n{\omega}\lp^{p}+\n{d\omega}\lp^{p}\right)^{1/p}.
\end{eqnarray*}
La cohomologie du complexe $(\Omega^{*,p}(M),d)$ s'appelle la {\em cohomologie} $L^p$ de $M$. Elle est int{\'e}ressante surtout si $M$ est non compacte. 

\ali

Par d\'efinition, la cohomologie $L^p$ est invariante par diff\'eomorphisme bilipschitzien. Dans la classe des vari\'et\'es simplement connexes \`a courbure n\'egative ou nulle, c'est un invariant de quasiisom\'etrie (cf. \cite{Gromovasympt}). 

En toute g{\'e}n{\'e}ralit{\'e}, la cohomologie $L^p$ se d{\'e}compose en
cohomologie r{\'e}duite et torsion
$$ 
0\to \tp{*}\to\hp{*}\to\rp{*}\to 0
$$
o{\`u} la {\em cohomologie r{\'e}duite} est $\rp{*}=\ker d/\overline{\im d}$ et la {\em torsion} est $\tp{*}=\overline{\im d}/\im d$. La cohomologie
r{\'e}duite (parfois notée $\overline{H}^{k}_{(p)}$) est un espace de Banach sur lequel les isom{\'e}tries de $M$ agissent isom{\'e}triquement. La torsion est non s{\'e}par{\'e}e.

Par exemple, la cohomologie $L^p$ de la droite r\'eelle est enti\`erement de torsion. La cohomologie $L^p$ du plan hyperbolique est enti\`erement r\'eduite, voir au paragraphe \ref{torcontr}. N\'eanmoins, cohomologie r\'eduite et torsion coexistent souvent.

\subsection{Torsion et pincement}

La torsion en cohomologie $L^p$ est sensible au pincement de la courbure sectionnelle, de façon optimale.

\label{A}
\begin{stheo} {\em \cite{Ppince}.}
Soient $\delta\in]-1,0[$ un r{\'e}el, $n$ et $k=2,\ldots,n$ des entiers. Soit $M$ une vari{\'e}t{\'e} riemannienne compl{\`e}te de dimension $n$, simplement connexe, dont la courbure sectionnelle $K$ satisfait $-1\leq K\leq\delta$. Alors
$$
\tp{k}(M)=0, \quad\hbox{i.e.}\quad \hp{k}(M)\quad\hbox{est s{\'e}par{\'e} pour}\quad 1<p<1+\frac{n-k}{k-1}\sqrt{-\delta}.
$$
\end{stheo}

Des exemples homogènes montrent que la borne du Théorème \ref{A} est optimale pour toutes les valeurs de $n$, $k$, $\delta$. Les espaces symétriques de rang un à courbure non constante sont de bons candidats lorsque $k$ est la dimension du corps de base, $2$, $4$ ou $8$. Toutefois, il s'avère que leur cohomologie $L^p$ reste séparée au-delà de l'intervalle prescrit par le théorème \ref{A}, pour $k=2$, $4$ ou $8$ et $\delta=-\frac{1}{4}$.

\subsection{Torsion des espaces sym{\'e}triques de rang un}

\begin{theo}
\label{torsion}
Soit $M$ un espace sym{\'e}trique de rang un {\`a} courbure non constante, de dimension $n=m\f$, $\f=2$, $4$ ou $8$. On voit $M$ comme un produit semi-direct $\R\ltimes_{\alpha}N$, où $\alpha$ a deux valeurs propres, 1 et 2. Etant donn{\'e} $k\leq \tr\alpha/2$, notons ${\rm suiv}\,\sigma(k)$ le second {\'e}l{\'e}ment, dans l'ordre décroissant, de l'ensemble $\sigma(k)$ des valeurs propres de $\L^{k}\alpha^{\top}$. Si $\tr\alpha/p >{\rm suiv}\,\sigma(k-1)$, alors $\tp{k}(M)=0$ sauf peut-{\^e}tre lorsque $\tr\alpha/p=\max\sigma(k-1)$.
\end{theo}

Par exemple, si $k=\f$, l'intervalle du théorème \ref{A} est $[1,\tr\alpha/\max\sigma(k-1)[$. La torsion $L^p$ reste donc nulle au-delà de cette borne. Il reste à traiter le cas limite $p=\tr\alpha/\max\sigma(\f-1)=(m\f +\f-2)/(2\f-2)$. Il est vraisemblable que la torsion est nulle pour cette valeur particulière, inaccessible par notre méthode. C'est le cas pour le plan hyperbolique complexe ($\f=2$, $m=2$). En effet, il s'agit des 2-formes $L^2$, cas où diverses méthodes sont disponibles, voir \cite{Borel}, \cite{Gromovkahler}. Mieux, N. Lohoué, a prouvé que la cohomologie $L^p$ reste séparée pour $p$ voisin de 2, lorsque $0$ est isolé dans le spectre du laplacien $L^2$, d'abord dans le cas des espaces symétriques, \cite{Lohoue98}, puis dans le cas général, \cite{Lohoue08}. Son résultat recouvre donc le nôtre dans ce cas particulier. 

\subsection{Cup-produit en cohomologie \texorpdfstring{$L^p$}{}}

Le produit extérieur d'une $k$-forme fermée $L^q$ et d'une $\ell$-forme fermée $L^r$ est une $(k+\ell)$-forme fermée $L^p$ si $\frac{1}{p}=\frac{1}{q}+\frac{1}{r}$. Si l'une des deux est exacte, il en est de même du produit. Le cup-produit $\smile$ est donc bien défini en cohomologie $L^p$ : il va de $H^{k,q}(M)\times H^{\ell,r}(M)$ dans $H^{k+\ell,p}(M)$. 

Cette opération est clairement fonctorielle sous les difféomorphismes bilipschitziens. Nous aurons besoin d'un résultat un peu plus général.

\label{qi}
\begin{stheo}
{\em \cite{qi}}
Soient $M$, $M'$ des variétés riemanniennes complètes, simplement connexes, à courbure pincée négativement. Soit $f:M\to M'$ une quasiisométrie. Alors $f$ induit un isomorphisme $f^{*}:\hp{*}(M')\to\hp{*}(M)$ qui préserve le cup-produit.
\end{stheo}

\subsection{Algèbre de Royden}

\begin{defi}
{\em (Bourdon-Pajot \cite{BP}, adapté au cas des variétés}). Soit $M$ une variété riemannienne. Soit $q\geq 1$. L'{\em algèbre de Royden} $\mathcal{R}^q (M)$ est l'espace des fonctions bornées $u$ sur $M$ telles que $du\in L^q$.
\end{defi}

Par définition, tout $u\in \mathcal{R}^q (M)$ définit une classe notée
$[du]\in H^{1,q}(M)$. Si $M$ est simplement connexe, l'application
$u\mapsto[du]$, $\mathcal{R}^q (M)\to H^{1,q}(M)$ est surjective, et son
noyau est $\R\oplus I$ où $I$ désigne l'idéal des fonctions $L^q \cap
L^{\infty}$. Ce lien avec la cohomologie $L^q$ rend les algèbres de Royden
fonctorielles sous les quasiisométries, voir \cite{BP}.

\begin{theo}
\label{subalgebra}
Soit $M$ une vari{\'e}t{\'e} riemannienne compl{\`e}te de dimension $n$,
simplement connexe, dont la courbure sectionnelle $K$ satisfait $-1\leq
K\leq\delta<0$. Soient $r>p\geq 1$, soit $q$ le réel tel que $\frac{1}{p}=\frac{1}{q}+\frac{1}{r}$. Supposons que
\begin{eqnarray*}
p<1+\frac{n-k}{k-1}\sqrt{-\delta}\quad\mathrm{et}\quad r<1+\frac{n-k+1}{k-2}\sqrt{-\delta}.
\end{eqnarray*}
Soit $\kappa\in H^{k-1,r}(V)$. Alors l'ensemble 
\begin{eqnarray*}
\mathcal{R}_{\kappa}=\{u\in\mathcal{R}^q (M)\,|\,[du]\smile\kappa=0\}
\end{eqnarray*}
est une sous-algèbre de $\mathcal{R}^q (M)$.
\end{theo}

\begin{rem}
Soit $\delta<-\frac{1}{4}$, $\f=2$, $4$ ou $8$, $n=m\f+\f-2$, $k=\f$, $q=r=2p$. Il existe $p_0 (\delta)>(m\f +\f-2)/(2\f-2)$ tel que le Théorème \ref{subalgebra} s'applique pour tout $p<p_0$.
\end{rem}

La borne sur $p$ donnée dans le théorème \ref{subalgebra} est optimale,
notamment pour le plan hyperbolique complexe. Plus précisément, on doit
utiliser une version localisée dans certains ouverts, voir la Proposition
\ref{preuvesubalgebra} pour un énoncé précis. De même, on a simplifié
l'énoncé suivant, dont la version technique est la Proposition \ref{nosub}.

\begin{theo}
\label{nosubalgebra}
Soit $M=\mathbb{C}H^2$ le plan hyperbolique complexe. Soit $p$ tel que
$2<p<4$. Soit $q=r=2p$. Soit $\eta$ un point de la sphère à l'infini. Il existe une forme fermée $\omega$ et une fonction $u$ sur $M$ telles que, pour tout cône $V$ de sommet $\eta$ sur un ouvert de la sphère à l'infini dont l'adhérence ne contient pas $\eta$,
\begin{enumerate}
  \item $\omega\in \Omega^{1,p}(V)$ et $u\in\mathcal{R}^{q}(V)$ ;
  \item $[du]\smile[\omega]=0$ dans $\hp{2}(V)$ ;
  \item $d(u^2)\smile[\omega]\not=0$ dans $\hp{2}(V)$.
\end{enumerate} 
\end{theo}

C'est le théorème \ref{nosubalgebra} qui ne s'étend pas aux autres espaces
symétriques de rang un. Le mécanisme algébrique qui s'y oppose est expliqué
dans le Corollaire \ref{cordb}. 

\subsection{Cup-produit et pincement}

En combinant les (versions techniques des) théorèmes \ref{qi}, \ref{subalgebra} et \ref{nosubalgebra}, on obtient le pincement optimal pour le plan hyperbolique complexe.

\begin{cor}
\label{thmpince}
Soit $M$ une vari{\'e}t{\'e} riemannienne compl{\`e}te de dimension $n$,
simplement connexe, dont la courbure sectionnelle $K$ satisfait $-1\leq K\leq\delta<0$. Si $\delta<-\frac{1}{4}$, $M$ n'est pas quasiisométrique au plan hyperbolique complexe $\mathbb{C}H^2$.
\end{cor}

\section{Méthode}

\subsection{Formule de Künneth}

L'idée de \cite{Ppince} est qu'une forme fermée $L^p$ sur une variété riemannienne à courbure sectionnelle négative pincée possède une valeur au bord, qui est une forme fermée de même degré sur le bord à l'infini. Cette valeur au bord ne dépend que de la classe de cohomologie $L^p$. En degré 1, cette construction marche pour tout $p\geq 1$. En degré supérieur, il faut supposer que $p$ est assez petit. \cite{Ppince} donne une borne en fonction du pincement de la courbure, celle du Théorème \ref{A}. L'opérateur valeur au bord injecte l'espace $\hp{k}$ dans un espace fonctionnel de formes différentielles fermées sur le bord, cela entraîne que $\hp{k}$ est séparé. Par exemple, pour le plan hyperbolique complexe, $\hp{2}$ est séparé pour $p<2$.

Pour aller au-delà de la borne du Théorème \ref{A}, on développe un avatar de la formule de K{\"u}nneth. En effet, une fois fixé un point à l'infini, la variété devient un produit $\R\times H$ où $H$ est une horosphère. Toutefois, la métrique riemannienne sur ce produit n'est pas un produit riemannien, ni même un produit tordu (le cas des produits tordus est bien compris, voir \cite{GKS,Kop1,Kop2,Kop3,Kop4,KS}). Lorsque la variété est homogène, il s'agit d'un {\em produit pluri-tordu}, au sens où la métrique croît exponentiellement, mais avec des taux de croissance différents suivant les directions. Du coup, la formule à la Künneth n'exprime pas la cohomologie $L^p$ de $\R\times H$ en fonction de la cohomologie $L^p$ de chaque facteur, mais elle en ramène le calcul à celui de la cohomologie d'un complexe
$\bp{*}$ de formes diff{\'e}rentielles sur $H$. L'espace $\bp{*}$ est
constitu{\'e} de formes diff{\'e}rentielles dont certaines composantes s'annulent et les autres appartiennent {\`a} des sortes d'espaces de Besov
anisotropes. Par exemple, l'espace hyperbolique r{\'e}el ({\`a} courbure
constante) s'{\'e}crit $\R\times \R^{n-1}$ et $\bp{k}$ est {\`a} peu de choses
pr{\`e}s l'espace des $k$-formes diff{\'e}rentielles ferm{\'e}es sur $\R^{n-1}$ {\`a} coefficients dans l'espace de Besov $\d B_{p,p}^{-k+(\frac{n-1}{p})}$, voir \cite{Pweb}. Autrement dit, la question de savoir si la cohomologie $L^{p}$ est s{\'e}par{\'e}e ou m{\^e}me nulle m{\^e}le analyse et alg{\`e}bre. 

Lorsque $p=2$, la torsion est li{\'e}e {\`a} la pr{\'e}sence de spectre
proche de $0$ pour le laplacien. Le calcul de la
torsion pour des produits tordus s'apparente {\`a} l'{\'e}tude des
petites valeurs propres du laplacien sur les formes dans la limite
adiabatique, voir par exemple \cite{MazzeoMelrose} ou \cite{Forman}, o{\`u} les consid{\'e}rations alg{\'e}briques jouent un grand r{\^o}le. La nouveaut{\'e} tient dans le fait que les poids (la ``taille'' des fibres) tendent simultan{\'e}ment vers l'infini et vers $0$ suivant les directions (m{\'e}lange de limite adiabatique et antiadiabatique). C'est le caract{\`e}re antiadiabatique qui produit les espaces de Besov.

La preuve de la formule de Künneth est inspirée des travaux de A.N. Livsic \cite{liv} sur la cohomologie des systèmes dynamiques. Nous nous sommes aussi inspirés de la terminologie des systèmes dynamiques : l'hypothèse qui fait marcher la formule de Künneth pluri-tordue en degré $k$, pour un exposant $p$ est l'existence d'un flot $(k-1,p)$-Anosov. Un flot $(0,\infty)$-Anosov n'est autre qu'un flot d'Anosov.

\subsection{Torsion}

Les formes basiques sur $\R\times H$ sont les formes différentielles sur $H$ qui, quand on les tire en arrière sur $\R\times H$, possèdent (ainsi que leur différentielle extérieure), $-1$ dérivée dans $L^p$. Cette propriété d'allure analytique entraîne aussi l'annulation pure et simple de certaines composantes.

Dans le langage de la formule de Künneth pluri-tordue, voici ce qui se produit dans le régime du Théorème \ref{A} : le complexe des formes basiques $\bp{*}$ est de longueur 1, il n'est non nul qu'en un seul degré $k$ (dépendant de $p$). Par conséquent, sa cohomologie, trivialement séparée, est égale à l'espace $\bp{k}$, qui ne contient que des formes fermées. 

Pour les espaces hyperboliques complexes, le complexe des formes basiques est souvent de longueur 2. C'est le cas, pour le plan hyperbolique complexe, lorsque $2<p<4$. Les espaces non nuls sont $\bp{1}$ et $\bp{2}$. Ici, $H$ s'identifie au groupe d'Heisenberg de dimension 3, avec sa forme de contact invariante à gauche $\tau$. Les éléments de $\bp{1}$ sont des $1$-formes différentielles de la forme $\psi=u\tau$, où $u$ est une fonction appartenant à un espace de Besov anisotrope. Leurs différentielles s'écrivent $du\wedge\tau+u\,d\tau$, i.e. elles se décomposent en deux morceaux, $d_{0}\psi=u\,d\tau$ et $d_{1}\psi=du\wedge\tau$, le second s'exprimant au moyen des dérivées du premier. En voici deux conséquences.
\begin{enumerate}
  \item Cette propriété n'est pas automatiquement satisfaite par une $2$-forme fermée sur $H$. Autrement dit, une $2$-forme fermée de $\bp{2}$ est exacte si et seulement si elle satisfait une équation différentielle non triviale.
  \item $\psi$ est déterminée algébriquement par $d\psi$. Cela entraîne que l'image $d\bp{1}$ est fermée dans $\bp{2}$.
\end{enumerate}
Le point 2 se généralise aux autres espaces symétriques de rang un, dans le degré égal à la dimension ($2$, $4$ ou $8$) du corps de base (et aussi bien d'autres degrés, voir le Théorème \ref{torsion}). En revanche, le phénomène 1 ne se produit pas dans ce degré. Il se produit dans d'autres degrés, mais c'est sans utilité pour la question du pincement optimal.

Pour une autre approche de la torsion dans le cas des espaces hyperboliques complexes, voir \cite{Prumin}.

\subsection{Algèbre de Royden}

Dans le régime du Théorème \ref{A}, une forme fermée $L^p$ est exacte si et seulement si sa valeur au bord est nulle. Il s'agit d'une condition non différentielle. L'espace des solutions est donc un module sur une algèbre de fonctions sur le bord (l'algèbre de Royden, introduite par M. Bourdon et H. Pajot, fait l'affaire) ? Pas exactement, puisque les valeurs au bord sont automatiquement fermées. Mais nous parvenons tout de même à saisir le changement, lorsque $p$ dépasse la borne du Théorème \ref{A}, de la caractérisation des formes exactes, au moyen de la structure multiplicative sur la cohomologie $L^p$. Si $u$ appartient à la bonne algèbre de Royden, et $\omega$ est une forme fermée, $du\wedge\omega$ est exacte si et seulement si sa valeur au bord $du_{\infty}\wedge\omega_{\infty}$ est nulle. Les fonctions $u_{\infty}$ satisfaisant cette propriété forment une sous-algèbre, c'est le Théorème \ref{subalgebra}.

Dans le cas du plan hyperbolique complexe, lorsque $2<p<4$, $\omega$ est de degré 1, c'est la différentielle d'une fonction $v$. Alors $du\wedge dv$ est exacte si et seulement si, en (presque) tout point du bord à l'infini, $du\wedge dv$ s'annule sur le plan de contact. De ce système d'equations différentielles bilinéaire, du premier ordre, sur le groupe d'Heisenberg, nous n'avons pas su trouver de solutions qui décroissent à l'infini. Ne disposant que d'un contre-exemple local à la propriété de sous-algèbre, nous avons dû donner une version localisée (au voisinage d'un point à l'infini) du Théorème \ref{subalgebra}, c'est la proposition \ref{preuvesubalgebra}. Avec l'invariance de la cohomologie $L^p$ (localisée) par quasiisométrie, on obtient le pincement optimal du plan hyperbolique complexe.

\subsection{Plan de l'article}

La formule de K{\"u}nneth ne fonctionne pas pour certaines valeurs de $p$ baptis{\'e}es {\em exposants critiques}, et calcul{\'e}s dans des exemples en section \ref{critique}. La formule de K{\"u}nneth pluri-tordue est {\'e}nonc{\'e}e et démontrée en \ref{fkg}. Les sections \ref{lp+dlp} et \ref{res} tentent, de façon inparfaite, de cerner les espaces $\bp{k}$. En section \ref{tor}, on en tire un crit{\`e}re d'annulation de la torsion qui culmine avec la preuve du Th{\'e}or{\`e}me \ref{torsion} en \ref{preuvetorsion}. On trouve en section \ref{cup} les preuves des Théorèmes \ref{subalgebra} et \ref{nosubalgebra}, ou plutôt, de versions localisées de ces résultats. La preuve du Corollaire \ref{thmpince} s'en déduit aisément en section \ref{interqi}.

\section{Exposants critiques}
\label{critique}

On note $Gl^{+}(n-1)$ le groupe des matrices $n-1\times n-1$ de d{\'e}terminant strictement positif. On consid{\`e}re des repr{\'e}sentations
lin{\'e}aires $\rho$ du groupe $Gl^{+}(n-1)$ sur des espaces euclidiens, qui
ont la propri{\'e}t{\'e} que le sous-groupe $SO(n-1)$ agit par isom{\'e}tries. Si $M$ est une vari{\'e}t{\'e} orient{\'e}e de dimension $n$ portant un champ de vecteurs $\xi$ qui ne s'annule pas, on construit le fibr{\'e} vectoriel
$E_\rho$ associ{\'e} au fibr{\'e} des rep{\`e}res directs ayant $\xi$ comme premier vecteur. Tout diff{\'e}omorphisme pr{\'e}servant l'orientation et le champ $\xi$ agit sur le fibr{\'e} $E_\rho$. Si de plus $M$ porte une m{\'e}trique riemannienne, alors $E_{\rho}$ porte une m{\'e}trique lui aussi. On note $\L^k$ la représentation de $Gl^{+}(n-1)$ sur $\L^k\R^{n-1}$.

\begin{defi}
Soit $M$ une vari{\'e}t{\'e} riemannienne de dimension $n$. Soit $\xi$ un champ de vecteurs unitaire complet sur $M$. Soit $\rho$ une repr{\'e}sentation lin{\'e}aire du groupe $Gl^{+}(n-1)$. On dit que $\xi$ est {\em $\rho$-Anosov} s'il existe des constantes positives $C$ et $\eta$ et une d{\'e}composition
$$
E_{\rho}=E_{\rho}^{+}\oplus E_{\rho}^{-}
$$ 
invariante par le flot $\phi_{t}$ de $\xi$, et telle que 
\begin{itemize}
\item quasi-orthogonalit{\'e} : si $e_{\pm}\in E_{\rho}^{\pm}$, alors
$$
|e_{+}|^{2}+|e_{-}|^{2}\leq C\,|e_{+}+e_{-}|^{2} ;
$$
\item contraction : si $e\in E_{\rho}^{\pm}$, alors pour tout $t\geq 0$,

$$
|\phi_{\pm t}(e)|\leq C\,e^{-\eta t}.
$$
\end{itemize}

Soit $p\geq 1$. On dit que $\xi$ est {\em $(k,p)$-Anosov}, ou bien que $p$ est un
{\em exposant non critique en degr{\'e} $k$} si $\xi$ est $\rho$-Anosov pour
la repr{\'e}sentation $\rho=\Lambda^{k}\otimes(det)^{-1/p}$ de
$Gl^{+}(n-1)$.

On dit que $\xi$ est {\em $(k,p)$-contractant} si $\xi$ est $(k,p)$-Anosov et $\L^{k}_{+}=0$. On dit que $\xi$ est {\em $(k,p)$-dilatant} si $-\xi$ est $(k,p)$-contractant.
\end{defi}

Concr{\`e}tement, pour $\rho=\Lambda^{k}$, les sections du fibr{\'e} dual
$(E_\rho)^{*}$ s'identifient aux $k$-formes dif\-f{\'e}\-ren\-ti\-el\-les sur $M$ annul{\'e}es par le produit int{\'e}rieur par $\xi$, $\iota_\xi$. Une d{\'e}composition $(k,p)$-Anosov est une d{\'e}composition quasi-orthogonale
$$
\L^{k}T^{*}M\cap\ker\iota_{\xi}=\L^{k}_{+}\oplus\L^{k}_{-}
$$
telle que 
$$
|(\L^{k}d\phi_{\pm t})_{|\L_{\pm}}|^{p}\leq \con e^{-p\eta t}Jac(\phi_{t})
$$
pour tout $t\geq 0$. Si $\beta$ est une $k$-forme différentielle sur $M$ qui, en chaque point, est à valeurs dans $\L_{-}$, alors, pour tout $t\geq 0$,
\begin{eqnarray*}
\parallel \phi_{t}^{*}\beta\parallel_{L^{p}}\leq \con \,e^{-\eta t}\parallel \phi_{t}^{*}\beta\parallel_{L^{p}}.
\end{eqnarray*}

\remarque 
Un champ de vecteurs est Anosov au sens des syst{\`e}mes dynamiques
si et seulement si il est $(0,+\infty)$-Anosov.

\remarque Si la divergence $div(\xi)$ est born{\'e}e, alors les exposants non
critiques forment un ouvert de $]1,+\infty]$.

\exemple S'il existe $\eta>0$ tel que $div(\xi)\geq\eta$, alors aucun
exposant $p\geq 1$ n'est critique en degr{\'e}s $0$, seul 1 est critique en degré $n-1$. Le champ $\xi$ est $(0,p)$-contractant pour tout $p\geq 1$ et $(n-1,p)$-dilatant pour tout $p>1$.

En effet, si $k=0$ ou $k=n-1$, $\L^{k}T^{*}M\cap\ker\iota_{\xi}$ est de
dimension $1$, $\L^{0}d\phi_{t}$ est l'identit{\'e}, $|\L^{n-1}d\phi_{t}|=Jac(\phi_{t})\geq e^{\eta t}$.

Les deux paragraphes suivants donnent des exemples moins triviaux.

\subsection{Variétés à courbure négative pincée}

Dans une variété simplement connexe à courbure négative, un {\em champ de vecteurs de Busemann} est la limite (lorsqu'elle existe) des gradients des fonctions distance à une suite de points tendant vers l'infini. Un tel champ $\xi$ est aussi différentiable que la métrique. Son flot augmente exponentiellement la longueur des vecteurs, donc $\xi$ est $(k,+\infty)$-dilatant pour tout $k=1,\ldots,n$. On s'attend à ce qu'il reste $(k,p)$-dilatant pour $p$ assez grand. La proposition suivante le garantit, sous une hypothèse faisant intervenir les deux bornes de la courbure sectionnelle.

\begin{prop}
\label{pince}
{\em (Voir \cite{Ppince}, Proposition 5)}.
Soit $M$ une vari{\'e}t{\'e} riemannienne compl{\`e}te de dimension $n$,
simplement connexe, dont la courbure sectionnelle $K$ satisfait $-1\leq K\leq\delta<0$. Soit $\xi$ un champ de vecteurs de Busemann. Si $k=0,\cdots,n-1$ et si $p>1$ satisfait
$$
p<1+\frac{n-k-1}{k}\sqrt{-\delta},\quad\hbox{ (resp. }
p>1+\frac{n-k-1}{k\sqrt{-\delta}}),
$$
alors le champ $\xi$ est $(k,p)$-contractant (resp. $(k,p)$-dilatant).
\end{prop}

\subsection{Groupes de Lie}

\begin{prop}
Soit $G$ un groupe de Lie muni d'une m{\'e}trique riemannienne
invariante {\`a} gauche. Soit $\xi$ un champ de vecteur invariant {\`a} gauche sur $G$. On note $\xi^{\bot}\subset\G^*$ le sous-espace des formes linéaires qui s'annulent sur $\xi$. On note $w_1 ,\ldots, w_{\dim(\mathfrak{g})-1}$ les parties réelles des valeurs propres de $ad_{\xi}^{\top}$ sur $\xi^{\bot}$, comptées autant de fois que leurs multiplicités. On note $\d h={\rm tr}\,ad_{\xi}=\sum_{j=1}^{\dim(\mathfrak{g})-1} w_j$. Un r{\'e}el $p\geq 1$ est un exposant critique en degr{\'e} $k$ pour $\xi$ si et seulement si $h/p$ s'écrit comme somme de $k$ nombres parmi les $w_j$. 
\end{prop}

\preuve
Soit $\L^{k}_{+}\subset \L^{k}\G^{*}\cap\ker\iota_{\xi}=\L^{k}\xi^{\perp}$
(resp. $\L^{k}_{-}$) la somme des espaces caract{\'e}ristiques de
$\L^{k}ad_{\xi}^{\top}$ relatifs {\`a} des valeurs propres $w$ telles que $h-pw>0$ (resp.
$h-pw<0$). On note de la m{\^e}me fa\c con les sous-fibr{\'e}s obtenus en
translatant {\`a}
gauche ces sous-espaces. Notons $\eta_{0}$ le plus petit des nombres
$|h-pw|$ o{\`u} $w$
d{\'e}crit les valeurs propres de $\L^{k}ad_{\xi}$ sur $\L^{k}\xi^{\perp}$.
Fixons un
$\eta<\eta_{0}$.

Par d{\'e}finition,
$L_{\exp -t\xi}\circ R_{\exp t\xi}=Ad_{\exp t\xi}=\exp(t\,ad_{\xi})$. Comme les
sous-fibr{\'e}s
$\L^{*}_{\pm}$ sont invariants par $ad_{\xi}$, ils sont invariants par les
translations {\`a} droite $\phi_{t}=R_{\exp t\xi}$ qui constituent le flot de
$\xi$.

Le champ de vecteurs $\xi$ a une divergence constante {\'e}gale {\`a} $h$, donc le jacobien de son flot vaut exactement $e^{ht}$.

Soit $\omega\in \L^{k}_{+}$ une forme diff{\'e}rentielle invariante
{\`a} gauche qui est dans l'image de $\iota_{\xi}$. On a
$$
(\phi_{t})^{*}\omega=Ad_{\exp t\xi}^{*}\omega=(\exp t\,ad_{\xi})^{*}\omega.
$$
Comme $\eta<\eta_{0}$, il existe une constante ind{\'e}pendante de $\omega$ et de $t$ telle que pour tout $t\geq 0$,
$$
|(\exp t\,ad_{\xi})^{*}\omega|\leq\con e^{t(h-\eta)/p}|\omega|.
$$
Il vient
$$
|(\phi_{t})^{*}\omega|^{p}=\con e^{-\eta t}J(\phi_{t})|\omega|^{p}.
$$
On conclut que $p$ est non critique en degr{\'e} $k$.

Inversement, s'il existe un couple de valeurs propres con\-ju\-guées
$(\lambda,\overline{\lambda})$ de $(\L^{k}ad_{\xi})_{|\L^{k}\xi^{\perp}}$ de partie r{\'e}elle $h/p$, alors sur la partie r{\'e}elle de la somme des sous-espaces caract{\'e}ristiques correspondants, $\n{\exp t\,ad_{\xi}}$ est polynomial en $t$ donc ne d{\'e}cro{\^\i}t pas exponentiellement, et $p$ est critique en degr{\'e} $k$.

Enfin, les valeurs propres de l'action de $\xi$ sur $\L^{k}(\xi^{\bot})\subset (\mathfrak{g}^{*})^{\otimes k}$ sont les sommes de $k$ valeurs propres de $ad_{\xi}$ sur $\xi^{\bot}$.\qed

\begin{cor}
Si $G$ n'est pas unimodulaire, un champ de vecteurs invariant {\`a} gauche
g{\'e}n{\'e}rique sur
$G$ a au plus un nombre fini d'exposants critiques.
\end{cor}

\begin{cor}
Soit $H$ un groupe de Lie d'alg{\`e}bre de Lie $\H$. Soit $\alpha$ une
d{\'e}rivation de $\H$, semi-simple, à valeurs propres réelles. Autrement dit, il existe une base de $\H$ dans laquelle la matrice de $\alpha$ est diagonale, $diag(w_1 ,\ldots,w_{n-1})$. Soit $G=\R\times H$ muni de la multiplication
$$
(t,h)(t',h')=(t+t',e^{t\alpha}(h)h').
$$
Soit $\xi=\frac{\partial}{\partial t}$ le champ de vecteurs invariant {\`a} gauche correspondant au facteur $\R$. Alors $\xi^{\bot}=\H^*$. Supposons $h=\mathrm{tr}(\alpha)>0$. Alors, pour tout $k=0,\ldots n-1=\dim(H)$ et $p\geq 1$, 
\begin{eqnarray*}
\L^{k}_{+(p)}=\bigoplus_{\{(w_{j_1},\ldots,w_{j_k})\,|\,w_{j_1}+\cdots+w_{j_k}>h/p\}}\bigotimes_{\ell=1}^{k}\L_{w_{j_\ell}},
\end{eqnarray*}
où, par convention, lorsqu'un même nombre $w$ est répété $m$ fois, le produit tensoriel doit être interprété comme une puissance extérieure $m$-ème.
\end{cor}

\subsection{Exemples}

\exemple Soit $G=\R\ltimes_{\alpha}\R^{2}$ le produit semi-direct d{\'e}fini par la d{\'e}rivation $\alpha=\begin{pmatrix}\l&0\| 0&1\end{pmatrix}$.
Si $\l>0$, $G$ est {\`a} courbure sectionnelle strictement n{\'e}gative. $G$
est unimodulaire si et seulement si
$\l=-1$. Si $\l=0$, $G$ est le produit riemannien d'une droite et d'un plan
hyperbolique.
Noter que pour tout $\ell\not=0$, les groupes obtenus pour les valeurs
$\ell$ et $1/\ell$
de $\lambda$ sont isomorphes. Par cons{\'e}quent, sans perdre de
g{\'e}n{\'e}ralit{\'e}, on peut
se limiter aux valeurs de $\l\in[-1,1]$.

On s'int{\'e}resse au champ de vecteurs invariant {\`a} gauche
correspondant au facteur $\R$. Si $\l\not=-1$, il n'a d'exposants critiques
qu'en degr{\'e} 1,
il en a deux (ce sont $p=1+\inv{\l}$ et $p=1+\l$) si $\l\not=0,~1$, un (c'est
$p=2$) si $\l=1$ et
aucun si $\l=0$.

Si $\l=-1$, tout exposant $p\geq 1$ est critique en degr{\'e}s 0 et 2. En degr{\'e} 1, aucun exposant n'est critique.

\exemple Espace hyperbolique r{\'e}el. Soit $G$ le groupe
des dilatations et translations de $\R^{n-1}$. Alors $G$ agit simplement transitivement sur l'espace hyperbolique r{\'e}el de dimension $n$. Soit $\xi$ le champ de vecteurs invariant {\`a} gauche correspondant aux dilatations. Sur $\xi^{\perp}=\R^{n-1}$, $ad_{\xi}$ est l'identit{\'e}. Il y a exactement un exposant critique pour chaque degr{\'e} $k=1,\ldots,n-1$, c'est $p=\frac{n-1}{k}$. Il n'y en a pas en degr{\'e} $0$. Dans le tableau, on note $\L^k_k =\L^{k}(\R^{n-1})^*$.
\begin{center}
\begin{tabular}{|c|ccccc|}
\hline
 $p$  &$1$&&$(n-1)/k$&&$+\infty$\\\hline
$\L^{k}_{+(p)}$&&$0$&&$\L^{k}_{k}$&\\\hline
$\L^{k}_{-(p)}$&&$\L^{k}_{k}$&&$0$&\\\hline
\end{tabular}
\end{center}
Autrement dit, $\xi$ est $(k,p)$-contractant si et seulement si $k=0$ (pour tout $p$) ou $1\leq k\leq n-2$ et $p<\frac{n-1}{k}$. $\xi$ est $(k,p)$-dilatant si et seulement si $k\geq 1$ et $p>\frac{n-1}{k}$.

\exemple Espace hyperbolique complexe. Soit $G$ le groupe des dilatations
translations du groupe d'Heisenberg $N$ de dimension $2m-1$. Alors $G$
agit simplement transitivement sur l'espace hyperbolique complexe de dimension r{\'e}elle $n=2m$, i.e., la boule de $\C^m$ munie de sa m{\'e}trique de Bergmann. Soit $\xi$ le champ de vecteurs invariant {\`a}
gauche correspondant aux dilatations. Sur $\xi^{\perp}=\N$, $ad_{\xi}$ induit une graduation $\N=\N_1 \oplus \N_{2}$. $ad_{\xi}$ vaut
l'identit{\'e} sur $\N_1$ et $2$ fois l'identit{\'e} sur $\N_{2}$. Une graduation des formes différentielles invariantes à gauche en découle. On numérote les sous-espaces par la valeur propre correspondante de $\L^{k}ad_{\xi}^{\top}$. On a donc, pour $k=1,\ldots,2m-2$,
\begin{eqnarray*}
\L^k \N^*=\L^k_k \oplus \L^k_{k+1} =\L^k \N_{1}^* \oplus \L^{k-1}\N_{1}^* \otimes (\N_{2})^* .
\end{eqnarray*}
et, pour les valeurs extrêmes,
\begin{eqnarray*}
\L^{0}\N^*=\L^{0}_{0}=\R,\quad \L^{2m-1}\N^* =\L^{2m-1}_{2m}=\L^{2m-2}\N_{1}^* \otimes (\N_{2})^* .
\end{eqnarray*}
Aucun exposant n'est critique en degr{\'e} $0$. Seul $1$ est critique en degr{\'e} $2m-1$. Par conséquent, $\L^{0}_{+(p)}=0$, $\L^{0}_{-(p)}=\R$ pour tout $p\geq 1$, $\L^{2m-1}_{+(p)}=\L^{2m-1}_{2m}$, $\L^{2m-1}_{-(p)}=0$ pour tout $p> 1$. En chaque degr{\'e} $k=1,\ldots,2m-2$, il y a exactement deux exposants critiques, $p=2m/(k+1)$ et $p=2m/k$. D'où le tableau.
\begin{center}
\begin{tabular}{|c|ccccccc|}
\hline
 $p$  &$1$&&$2m/(k+1)$&&$2m/k$&&$+\infty$\\\hline
$\L^{k}_{+(p)}$&&$0$&&$\L^{k}_{k+1}$&&$\L^{k}_{k}\oplus\L^{k}_{k+1}$&\\\hline
$\L^{k}_{-(p)}$&&$\L^{k}_{k}\oplus\L^{k}_{k+1}$&&$\L^{k}_{k}$&&$0$&\\\hline
\end{tabular}
\end{center}
Autrement dit, $\xi$ est $(k,p)$-contractant si et seulement si $k=0$ (pour tout $p$) ou $1\leq k\leq 2m-2$ et $p<\frac{2m}{k+1}$. $\xi$ est $(k,p)$-dilatant si et seulement si $k\geq 1$ et $p>\frac{2m}{k}$.

\exemple Espace hyperbolique quaternionien ou octonionien. Soit $N$ le radical unipotent du parabolique minimal de $Sp(m,1)$, $m\geq 2$ (resp. de $F_{4}^{-20}$). Son algèbre de Lie $\N$ possède une dérivation $\alpha$ ayant exactement deux valeurs propres $1$ et $2$. Autrement dit, $\N=\N_1 \oplus \N_{2}$ où $\dim(\N_1)=4m-4$ (resp. 8) et $\dim(\N_{2})=3$ (resp. 7). Alors $G=\R\ltimes_{\alpha}N$ agit simplement transitivement sur l'espace hyperbolique quaternionien de dimension r{\'e}elle $n=4m$ (resp. le plan hyperbolique octonionien). Soit $\xi$ le champ de vecteurs invariant {\`a}
gauche correspondant au facteur $\R$, de sorte que $ad_{\xi}=\alpha$. La graduation des formes différentielles invariantes à gauche comporte $4$ (resp. 8) sous-espaces en degrés $k=3,\ldots,4m-4$ (resp. 7 et 8), moins pour les valeurs extrêmes. Posons $\f =4$ pour les quaternions et $\f =8$ (et $m=2$) pour les octonions. Alors, pour $k=\f -1,\ldots, \f m-\f $,
\begin{eqnarray*}
\L^k \N^*=\bigoplus_{j=0}^{\f -1}\L^k_{k+j} =\bigoplus_{j=0}^{\f -1} \L^{k-j}\N_{1}^* \otimes \L^j (\N_{2})^* .
\end{eqnarray*}
et, pour les valeurs extrêmes, 
\begin{eqnarray*}
\L^{k}\N^*=\bigoplus_{j=0}^{k}\L^k_{k+j} \quad \textrm{ si }k\leq \f -2,\quad
\L^{k}\N^*=\bigoplus_{j=k-\f m+\f }^{\f -1}\L^k_{k+j} \quad \textrm{ si }k\geq \f m-\f +1.
\end{eqnarray*}
$h=\mathrm{tr}(\alpha)=\f m+\f -2$. Les exposants critiques sont les $(\f m+\f -2)/(k+j)$, $j=\max\{0,k-\f m+\f \},\ldots,\min\{k,\f -1\}$. Il y en a $\f $ si $k=\f -1,\ldots, \f m-\f $, moins pour les valeurs extrêmes de $k$. On trouve que $\xi$ est $(k,p)$-contractant si et seulement si $k=0$ (pour tout $p$) ou $1\leq k\leq \f m-\f $ et $p<\frac{\f m+\f -2}{k+\min\{k,\f -1\}}$. $\xi$ est $(k,p)$-dilatant si et seulement si $k\geq 1$ et $p>\frac{\f m+\f -2}{k+\max\{0,k-\f m+\f \}}$. Le spectre de $\L^k \alpha$ est l'intervalle d'entiers $[\max\{k,2k-\f m+\f \},\ldots,\min\{2k,k+\f -1\}]$. Si $1\leq k\leq \f m-1$, son plus grand élément est $\max\sigma(k)=\min\{2k,k-1+\f \}$. Si $2\leq k\leq \f m-2$, le suivant, par ordre décroissant, est $\mathrm{suiv}\,\sigma(k)=\min\{2k+1,k-1+\f \}-1$. Ces formules sont aussi valables pour les espaces hyperboliques réels (avec $\f =1$) et complexes (avec $\f =2$).

\remarque
Sachant que les espaces symétriques de rang un à courbure non constante sont $-\frac{1}{4}$-pincés, on constate que l'intervalle $[1,1+\frac{n-k}{k-1}\sqrt{-\delta}[$ pour lequel la Proposition \ref{pince} garantit que les champs de vecteurs de Busemann sont $(k-1,p)$-contractants coïncide avec l'intervalle observé pour ces espaces, $[1,\frac{\f m+\f -2}{k+\min\{k,\f -1\}}[$, seulement lorsque $k=\f$.

De même, l'intervalle de $(k-1,p)$-dilatation de la Proposition \ref{pince} est optimal pour les espaces symétriques de rang un à courbure non constante seulement pour $k=n-\f+1$.

\section{Formule de K{\"u}nneth pluri-tordue}
\label{fkg}

Dans cette section, on d{\'e}montre le principal r{\'e}sultat technique de cet
article, le th{\'e}or{\`e}me \ref{thm}.

\subsection{L'op{\'e}rateur \texorpdfstring{$B$}{}}

\begin{prop}
\label{B}
Soit $p\geq 1$. Soit $M$ une vari{\'e}t{\'e} riemannienne, soit $\xi$ un champ de vecteurs unitaire complet et $(k,p)$-Anosov sur $M$, de flot $\phi_t$. On note $\L_{+}\oplus\L_{-}$ la d{\'e}composition $(k,p)$-Anosov des $k$-formes annul{\'e}es par $\iota_{\xi}$. Si $\omega\in\L^{k}T^{*}_{x}M$, on note $\omega_{\pm}$ la projection orthogonale de $\omega$ sur $\L_{\pm}T^{*}_{x}M$.
Alors les op{\'e}rateurs
$$
\overline{B}:\omega\mapsto\int_{-\infty}^{0} (\phi_{t})^{*}\omega_{-}\,dt -
\int_{0}^{+\infty}
(\phi_{t})^{*}\omega_{+}\,dt,
$$
$$
i=\iota_{\xi}=\hbox{ produit int{\'e}rieur par }\xi
$$
et
$$
B=\overline{B}i
$$
sont born{\'e}s sur les $k+1$-formes $L^p$. De plus, il satisfont $Bi=iB=0$, $i(1-dB)=0$.
\end{prop}
\preuve
Par hypoth{\`e}se, pour tout $x\in M$ et tout $t\geq 0$,
$$
|(\phi_{t})^{*}i\omega_{+}|^{p}(x)\leq\con e^{-\eta
t}|i\omega_{+}|^{p}(\phi_{t}x)
J_{\phi_{t}}(x)
$$
donc
\begin{eqnarray*}
\int_{M} |(\phi_{t})^{*}i\omega_{+}|^{p}(x)\,dx
&\leq&\con e^{-\eta t}
\int_{M} |\omega|^{p}(\phi_{t}x)J_{\phi_{t}}(x)\,dx\\
&=&\con e^{-\eta t}\int_{M} |\omega|^{p}(z)\,dz.
\end{eqnarray*}
Par cons{\'e}quent, la fonction $t\mapsto\n{(\phi_{t})^{*}i\omega_{+}}_{L^{p}}$ est
int{\'e}grable sur $\R_{+}$, la seconde int{\'e}grale converge dans
$L^{p}\Omega^{k}(M)$ et
$$
\n{(\phi_{t})^{*}i\omega_{+}}_{L^{p}}\leq \frac{1}{\eta}\n{i\omega_{+}}_{L^{p}}.
$$
Il vient $\d \n{B\omega}\leq \eta^{-1}\n{\omega}$.

L'op{\'e}rateur $i$ est {\'e}videmment born{\'e} sur $L^p$. Comme il commute avec
$\phi_t$, il
commute avec $\overline{B}$ d'o{\`u} $iB=i\overline{B}i=\overline{B}ii=0$. 

Notons ${\cal L}_{\xi}$ la
d{\'e}riv{\'e}e de Lie suivant $\xi$. Si $\omega$ est une $k+1$-forme lisse {\`a} support compact, alors pour tout $t$,
$$
\ddt(\phi_{t})^{*}\omega=(\phi_{t})^{*}{\cal L}_{\xi}\omega
={\cal L}_{\xi}(\phi_{t}^{*}\omega).
$$
D'autre part, au voisinage d'un point, $(\phi_{t})^{*}\omega=0$ pour $t$ assez grand. Par cons{\'e}quent
$$
{\cal L}_{\xi}\overline{B}\omega=\omega_{+}+\omega_{-}=\omega.
$$
Il vient 
$$
i\omega=i{\cal L}_{\xi}\overline{B}\omega={\cal L}_{\xi}B\omega
=(di+id)B\omega=idB\omega.
$$
Par densit{\'e}, cette identit{\'e} s'{\'e}tend {\`a} $L^{p}$.\qed

\subsection{Le complexe des formes basiques}

L'opérateur $B$ est notre candidat pour réaliser une homotopie de l'identité à un projecteur $P$ sur le sous-espace des formes annulées par $\iota_{\xi}$ et invariantes par le flot de $\xi$. Autrement dit, $1-P=dB+Bd$. Une difficulté surgit : $dB$ n'est pas borné sur $\op{*}(M)$. $P$ n'envoie pas $\op{*}(M)$ dans lui-même, mais dans un espace plus grand, l'espace $L^p +dL^p$.

\begin{defi}
\label{defpsi}
On note $\psip{k}(M)$ l'espace des $k$-formes différentielles $\psi$ telles qu'il existe $\beta\in\lpo{k}(M)$ et $\gamma\in\lpo{k-1}(M)$ telles que $\omega=\beta+d\gamma$. Il est muni de la norme
\begin{eqnarray*}
\n{\psi}_{\psip{*}}=\inf\{\left(\n{\beta}\lp^{p}+\n{\gamma}\lp^{p}\right)^{1/p}\,|\,\psi=\beta+d\gamma\}.
\end{eqnarray*}
\end{defi}

Le complexe $\psip{*}(M)$ calcule lui aussi la cohomologie $L^p$. En effet, si une forme $\omega\in\op{k}(M)$ s'écrit $\omega=d\psi$ où $\psi\in\psip{k-1}(M)$, i.e. $\psi=\beta+d\gamma$, alors $\beta\in\op{k-1}(M)$ et $\omega=d\beta$. Inversement, toute forme fermée $\psi\in\psip{k}(M)$, $\psi=\beta+d\gamma$, est cohomologue (au sens du complexe $\psip{*}$) à la forme fermée $\beta$, et $\beta\in\op{k}(M)$. Mais cela ne fournit pas directement une homotopie entre les deux complexes. Une telle homotopie sera construite plus loin, sous une hypothèse supplémentaire sur $M$, voir Proposition \ref{regularisation}.

\begin{defi}
\label{basique}
Soit $M$ une vari{\'e}t{\'e} riemannienne, $\xi$ un champ de vecteurs unitaire complet sur $M$. On note $\bp{*}(M,\xi)\subset \psip{*}(M)$ {\em le sous-complexe des formes qui sont annul\'ees par}
$\iota_{\xi}$ et $\iota_{\xi}\circ d$.
\end{defi}

Soit $G=\R\ltimes_{\alpha} H$ un produit semi-direct. Notons $\pi:G\to H$ la projection le long des orbites du champ de vecteurs $\xi$ correspondant {\`a} l'action {\`a} droite du facteur $\R$. L'op{\'e}ration $\pi^{*}$, image inverse par $\pi$, identifie $\bp{*}(G,\xi)$ {\`a} un sous-complexe de formes diff{\'e}rentielles sur $H$, qu'on note $\bp{*}$. La nature de cet espace ne saute pas aux yeux. Il s'av{\`e}re qu'il est souvent nul. Lorsqu'il ne l'est pas, c'est un espace fonctionnel constitu{\'e} de formes dont certains coefficients s'annulent (voir section \ref{res}) et les autres ont des d{\'e}riv{\'e}es fractionnaires anisotropes dans $L^p$ (voir section \ref{lp+dlp}).

\subsection{Théorème principal}

\begin{theo}
\label{thm}
Soit $M$ une vari{\'e}t{\'e} riemannienne. Soit $\xi$ un champ de vecteurs unitaire complet sur $M$. Soit $p\geq 1$, $1\leq k\leq n$. 
\begin{enumerate}
  \item Si le champ $\xi$ est $(k-1,p)$-Anosov, alors l'injection $I:\bp{k}(M,\xi)\subset \psip{k}(M)$ induit une surjection en cohomologie. 
  \item Si $k=1$ ou si $k\geq 2$ et $\xi$ est $(k-2,p)$-Anosov, alors $I$ induit une injection en cohomologie.
  \item Supposons de plus que $M$ est à géométrie bornée, que $p>1$ et que $\xi$ est $(j,p)$-Anosov pour tout $j=0,\cdots, n-1$. Alors les complexes $\bp{*}(M,\xi)$ et $\op{*}(M)$ sont homotopes.
\end{enumerate}
\end{theo}

\preuve
On utilise l'opérateur $B$ défini dans la Proposition \ref{B}.

Supposons que $p$ est non critique en degr{\'e} $k-1$. Alors $B$ est d{\'e}fini sur les $k$-formes $L^p$. Soit $\omega=\beta+d\gamma\in\psip{k}(M)$ une forme fermée. Alors $\delta=(1-dB)\beta$ est bien définie, et elle appartient à $\bp{k}(M,\xi)$. En effet, $\iota_{\xi}\delta=\iota_{\xi}(1-dB)\beta=0$ et $d\delta=d\beta=d\omega=0$. Comme $\delta=\beta-dB\beta=\omega-d(\gamma+B\beta)$, $\delta$ est cohomologue à $\omega$ dans $\psip{*}(M)$, donc $I$ est surjective en cohomologie. 

Supposons que $p$ est non critique en degr{\'e} $k-2$. Soit $\omega\in\bp{k}(M,\xi)$ une forme fermée. Supposons $\omega$ cohomologue à 0 dans $\psip{*}(M)$, $\omega=d(\beta+d\gamma)=d\beta$, où $\beta$ est $L^p$. Alors $\delta=(1-dB)\beta$ est bien définie, $\iota_{\xi}\delta=0$ et $d\delta=d\beta=\omega$ satisfait aussi $\iota_{\xi}d\delta=0$. Autrement dit, $\delta\in\bp{k-1}(M,\xi)$, $d\delta=\omega$ donc $\omega$ est cohomologue à 0 dans $\bp{*}(M,\xi)$, cela prouve que $I$ est injective en cohomologie.

Supposons désormais que $M$ est à géométrie bornée, et que $\xi$ est $(j,p)$-Anosov pour tout $j=0,\cdots, n-1$. On combine l'opérateur $P:\op{*}(M)\to\bp{*}(M,\xi)$ défini par $P=1-dB-Bd$ et l'opérateur de régularisation $S:\psip{*}(M)\to\op{*}(M)$ fourni par la proposition \ref{regul}. $P$ est borné de $\op{*}(M)$ dans $\psip{*}(M)$. Par définition, $Pd=dP$. En outre, comme $\iota_{\xi}B=0$, $\iota_{\xi}(1-dB)=0$, $\iota_{\xi}P=0$ et $\iota_{\xi}dP=\iota_{\xi}Pd=0$, donc $P$ est borné de $\op{*}(M)$ dans $\bp{*}(M,\xi)$. De son côté, $S$ s'écrit $S=1-dT-Td$ où $T$ est borné de $\omp{*}(M)$ dans $\op{*}(M)$. Par conséquent, $S$ est borné de $\bp{*}(M,\xi)$ dans $\op{*}(M)$, et
\begin{eqnarray*}
1-SP&=&1-S+S(dB+Bd)\\
&=&dT+Td+dSB+SBd\\
&=&d(T+SB)+(T+SB)d,
\end{eqnarray*}
où $T+SB$ est borné de $L^p$ vers $L^p$. De plus, $d(T+SB)=1-SP-(T+SB)d$ est borné de $\op{*}(M)$ vers $L^p$, donc $T+SB$ est borné de 
$\op{*}(M)$ dans lui-même. De même, 
\begin{eqnarray*}
1-PS&=&1-P+P(dT+Td)\\
&=&dB+Bd+dPT+PTd\\
&=&d(B+PT)+(B+PT)d.
\end{eqnarray*}
On constate que $\iota_{\xi}(B+PT)=0$, $\iota_{\xi}d(B+PT)=\iota_{\xi}(1-PS-(B+PT)d)=\iota_{\xi}=0$ sur $\ker \iota_{\xi}\cap\ker \iota_{\xi}d$. De plus
\begin{eqnarray*}
B+PT&=&B+(1-dB-Bd)T\\
&=&B+T-dBT-B(1-S-Td)\\
&=&T-dBT+BS+BTd
\end{eqnarray*}
est born{\'e} de $\psip{*}(M)$ dans lui-même, donc de $\bp{*}(M,\xi)$ dans lui-m{\^e}me. Les opérateurs $P$ et $S$ constituent donc une équivalence d'homotopie de complexes. \qed

\remarque
Comme la torsion $\tp{k}(M)$ est l'adhérence de 0 dans $\hp{k}(M)$, et $\rp{k}(M)=\hp{k}(M)/\tp{k}(M)$, ces deux espaces sont conservés par les isomorphismes. Sous les hypothèses 1 ou 2 du Théorème \ref{thm}, $\rp{k}(M)=R^k (\bp{*}(M,\xi))$ et $\tp{k}(M)=T^k (\bp{*}(M,\xi))$.

\subsection{L'opérateur \texorpdfstring{$P$}{}}

On n'utilisera pas directement l'énoncé 3 du théorème \ref{thm} dans cet article. Néanmoins, l'opérateur $P:\op{*}(M)\to\bp{*}(M,\xi)$ qui réalise l'homotopie est intéressant en lui-même. Il généralise la valeur au bord utilisée dans \cite{Ppince}. Ses propriétés, relatives au cup-produit notamment, seront détaillées en section \ref{cup}.

\begin{cor}
\label{corthm}
Soit $M$ une vari{\'e}t{\'e} riemannienne, non nécessairement complète. Soit $\xi$ un champ de vecteurs unitaire complet sur $M$. Soit $p\geq 1$, $1\leq k\leq n$. On suppose que $\xi$ est $(k-2,p)$ et $(k-1,p)$-Anosov. Alors les opérateurs 
\begin{eqnarray*}
P=1-dB-Bd :\op{k-1}(M)\to \bp{k-1}(M,\xi)
\end{eqnarray*}
et
\begin{eqnarray*}
P=1-dB:\op{k}(M)\cap\ker(d)\to \bp{k}(M,\xi)\cap\ker(d)
\end{eqnarray*}
sont bornés. Sur $\op{k-1}(M)$, $Pd=dP$. Enfin, $P$ induit un isomorphisme \break $\hp{k}(M)\to H^{k}(\bp{*}(M,\xi))$.
\end{cor}

\preuve
Elle combine la preuve du Théorème \ref{thm} et le fait que les complexes $\op{*}$ et $\psip{*}$ ont la même cohomologie.

Soit $\omega\in \op{k}(M)\cap\ker(d)$. Si $[P\omega]=0$ dans $H^{k}(\bp{*}(M,\xi))$, $P\omega=d(\beta+d\gamma)$ où $\beta$ et $\gamma$ sont $L^p$ et $\beta+d\gamma\in \bp{k-1}(M,\xi)$, alors $(1-dB)\omega=d(\beta+d\gamma)$, d'où
\begin{eqnarray*}
\omega=d(B\omega+\beta),
\end{eqnarray*}
et $B\omega+\beta\in\op{k-1}(M)$, donc $[\omega]=0$ dans $\hp{k}(M)$.

Soit $\psi=\beta+d\gamma\in \bp{k-1}(M,\xi)\cap\ker(d)$. Alors
\begin{eqnarray*}
\psi=P\beta+d(B\beta+\gamma-dB\gamma).
\end{eqnarray*}
Clairement, $B\beta+\gamma-dB\gamma\in\psip{k-1}(M)$, $\iota_{\xi}(B\beta+(1-dB)\gamma)=0$ et $\iota_{\xi}d(B\beta+(1-dB)\gamma)=\iota_{\xi}(\psi)=0$, donc $B\beta+\gamma-dB\gamma\in\bp{k-1}(M,\xi)$, et $[\psi]=[P\beta]$ in $H^{k}(\bp{*}(M,\xi))$.\qed

\medskip

A titre d'illustration, on prouve le résultat principal (Proposition 10) de \cite{Ppince}.

\begin{lemme}
\label{voisinage}
Soit $M$ une variété riemannienne complète, difféomorphe à $\R\times H$. Soit $U\subset H$ un ouvert, soit $V=\R\times U$ et $V_T =]T,+\infty[\times U$ (éventuellement, $T=-\infty$). Soient $p\geq 1$ et $k\geq 2$. Supposons que $\xi=\ddt$ est $(k-2,p)$-contractant et $(k-1,p)$-contractant. Alors l'opérateur $P$ est défini sur $\op{k-1}(V_T)$ et sur $\op{k}(V_T)\cap\ker(d)$. Il satisfait $Pd=dP$ sur $\op{k-1}(V_T)$ et donne un isomorphisme $\hp{k}(V_T)\to H^k (\bp{*}(V,\xi))$. Dans ce cas particulier, 
\begin{eqnarray*}
P\omega=\lim_{t\to+\infty}\phi_{t}^{*}\omega.
\end{eqnarray*}
\end{lemme}

\preuve
Comme $\xi$ est $(k-2,p)$-contractant et $(k-1,p)$-contractant, $\L^{k-2}_{-(p)}=0$ et $\L^{k-1}_{-(p)}=0$. Par conséquent, en degrés $k-1$ et $k$, la définition de l'opérateur $B$,
\begin{eqnarray*}
B\omega=-\int_{0}^{+\infty}(\phi_{t})^{*}(\iota_{\xi}\omega)_{+}\,dt
\end{eqnarray*} 
se généralise aux ouverts invariants par le semi-groupe $(\phi_t)_{t\geq 0}$. Le Corollaire \ref{corthm} se généralise aussi. 

$B$ est la limite quand $T\to+\infty$ des opérateurs
\begin{eqnarray*}
B_{T}\omega=-\int_{0}^{T} (\phi_{t})^{*}(\iota_{\xi}\omega)\,dt.
\end{eqnarray*}
Comme
\begin{eqnarray*}
\phi_{T}^{*}\omega-\omega=-dB_{T}\omega-B_{T}d\omega,
\end{eqnarray*}
il vient
\begin{eqnarray*}
P\omega=\omega-dB\omega-Bd\omega=\lim_{T\to+\infty}\phi_{T}^{*}\omega.\qed
\end{eqnarray*}

\subsection{Localisation}

Nous avons besoin d'une version du Corollaire \ref{corthm} localisée au
voisinage d'un point du bord à l'infini d'une variété à courbure
sectionnelle négative. On peut prendre de tels voisinages de la forme $\d
V_T =\bigcup_{t>T}\phi_{t}(U)=]T,+\infty[\times U$ où $U$ est un ouvert
relativement compact d'une horosphère. Le Théorème \ref{thm} (énoncés 1 et
2) s'applique à $V=\R\times U$, mais non à $V_T$. Néanmoins, le lemme suivant donne une information suffisante.

\begin{lemme}
\label{ouvert}
Soit $M$ une variété riemannienne complète, difféomorphe à $\R\times H$. Soit $U\subset H$ un ouvert borné contractile, à bord lisse, soit $V=\R\times U$ et $V_T =]T,+\infty[\times U$. Soient $p\geq 1$ et $k\geq 2$. Supposons que $\xi=\ddt$ est $(k-2,p)$- et $(k-1,p)$-Anosov. Soit $\omega\in\op{k}(V)$ une forme fermée qui est nulle sur $]-\infty,T+1[\times U$. Alors $[\omega_{|V_T }]\in \hp{k} (V_T )$ est nulle si et seulement si $[P\omega]\in H^k (\bp{*}(V,\xi))$ est nulle.
\end{lemme}

\preuve
Notons $Z=]T,T+1[\times U$. La suite exacte longue de la paire $(V_T ,Z)$ donne
\begin{eqnarray*}
\cdots\to \hp{k-1}(Z)\to \hp{k}(V_T ,Z)\to \hp{k}(V_T )\to \hp{k}(Z)\to\cdots
\end{eqnarray*}
Comme $Z$ est muni d'une métrique équivalente à la métrique produit sur $]T,T+1[\times U$, sa cohomologie $L^p$ coïncide avec la cohomologie ordinaire de $U$, qui est nulle, par hypothèse, en degrés $\geq 1$. Donc la flèche $j:\hp{k}(V_T ,Z)\to \hp{k}(V_T )$ est injective. Par hypothèse, $0=[\omega_{|V_T }]=j([\omega_{|(V_T ,Z)}])\in \hp{k}(V_T )$, donc $[\omega_{|(V_T ,Z)}]\in \hp{k}(V_T ,Z)$ est nulle, il existe $\psi\in\op{k-1}(V_T )$, nulle sur $Z$, telle que $d\psi=0$. On prolonge $\psi$ par 0. On obtient une forme $L^p$ sur $V$ telle que $d\psi=\omega$, i.e. $\omega$ est exacte. Réciproquement, si $\omega$ est exacte, sa restriction à $V_T$ l'est aussi. Il ne reste plus qu'à appliquer le Corollaire \ref{corthm} dans $V$.\qed

\section{Propriétés analytiques des formes basiques}
\label{lp+dlp}

Dans cette section, on décrit la norme des espaces $\psip{k}(M)$, en la reliant à une norme plus simple. Cela permettra, lors de la preuve du Théorème \ref{torsion}, de montrer qu'un opérateur candidat à être un inverse de la différentielle est borné sur $\bp{k}(M,\xi)$. Cette application reposera sur le Corollaire \ref{bornef}.

\subsection{L'espace \texorpdfstring{$\lmpo{k}(M)$}{}}
\label{lp-1}

La norme de l'espace $\psip{k}(M)$ est difficile {\`a} {\'e}valuer en
g{\'e}n{\'e}ral. On va construire une norme {\'e}quivalente plus commode, dans le cas particulier des vari{\'e}t{\'e}s riemanniennes {\`a} g{\'e}om{\'e}trie born{\'e}e.

\begin{defi}
Fixons un entier $\ell\geq 1$. Soit $M$ une vari{\'e}t{\'e} riemannienne. On dit que $M$ est {\`a} g{\'e}om{\'e}trie born{\'e}e s'il existe un $\epsilon>0$ tel que toutes les boules de rayon $\epsilon$ soient uniform{\'e}ment proches, au sens $C^\ell$, de la boule unit{\'e} de $\R^n$.
\end{defi}

\remarque Dans la litt{\'e}rature, on demande souvent seulement une borne
inf{\'e}rieure sur le rayon d'injectivit{\'e} et des bornes sur la courbure sectionnelle. Quitte {\`a} d{\'e}former la m{\'e}trique en une m{\'e}trique {\'e}quivalente plus r{\'e}guli{\`e}re (ce qui ne change pas la cohomologie $L^p$), on peut supposer que les d{\'e}riv{\'e}es covariantes jusqu'{\`a} l'ordre $\ell$ de la courbure sont uniform{\'e}ment born{\'e}es (voir \cite{BMR}). La proximit{\'e} $C^\ell$ ne co{\^u}te donc pas plus cher.

\begin{defi}
\label{locale}
{\em Soit $M$ une vari{\'e}t{\'e} riemannienne compacte (resp. compacte {\`a} bord). On note $\nabla$ la d{\'e}riv{\'e}e covariante. Soient $p$ et $p'$ des r{\'e}els positifs tels que $\d \inv{p}+\inv{p'}=1$. Soit $\omega$ une $k$-forme diff{\'e}rentielle sur $M$. Sa} norme $L_{-\ell}^{p}$ {\em est la borne sup{\'e}rieure des int{\'e}grales $\d\int_{M} \omega\wedge\phi$ sur les $n-k$-formes $\phi$ sur $M$, de classe $C^\ell$ (resp. nulles au voisinage du bord) et telles que}
$$
\n{\phi}_{L_{\ell}^{p'}(M)}:=(\n{\phi}^{p'}_{L^{p'}(M)}
+\n{\nabla\phi}^{p'}_{L^{p'}(M)}+\cdots+\n{\nabla^{\ell}\phi}^{p'}_{L^{p'}(M)})^{1/p'}\leq 1.
$$
\end{defi}

\remarque Le transport par un diff{\'e}omorphisme de classe $C^\ell$ est un
isomorphisme entre espaces $L_{-\ell}^{p}$.

\begin{defi}
\label{globale}
{\em Soit $M$ une vari{\'e}t{\'e} riemannienne {\`a} g{\'e}om{\'e}trie born{\'e}e. Soit $\omega$ une $k$-forme sur $M$. On} note
$$
\n{\omega}_{L_{-\ell}^{p}\Omega^{k}(M)}=(\int_{M}
\n{\omega_{|B(x,\epsilon)}}_{L^{p}_{-\ell}(B(x,\epsilon))}^{p}\,dx)^{1/p}.
$$
\end{defi}

\remarque Changer le rayon des boules conduit {\`a} une norme {\'e}\-qui\-va\-len\-te. Si $M$ est compacte, elle est automatiquement {\`a} g{\'e}om{\'e}trie born{\'e}e. Les d{\'e}finitions \ref{locale} et \ref{globale} donnent dans ce cas des normes {\'e}quivalentes.

\remarque Clairement, $L_{-\ell}^{p}\subset L_{-\ell-1}^{p}$ et la
diff{\'e}rentielle ext{\'e}rieure est born{\'e}e de $L_{-\ell}^{p}$ dans $L_{-\ell-1}^{p}$.

\begin{lemme}
\label{module}
Si $M$ est compacte ({\'e}ventuellement avec bord), l'espace $L_{-\ell}^{p}(M)$ est un mo\-du\-le sur l'algèbre $C^{\ell}(M)$.
\end{lemme}

\preuve
Soit $u$ une fonction de classe $C^\ell$ sur $M$, $\omega$ une $k$-forme
sur $B$ et $\phi$ une $n-k$-forme lisse {\`a} support compact dans $M$, telle que $\d\n{\phi}_{L_{\ell}^{p}(M)}\leq 1$. Alors
$\nabla u\phi=(\nabla u)\phi +u\nabla\phi$ entra{\^\i}ne
$$
\n{u\phi}_{L_{\ell}^{p'}}\leq \n{u}_{C^{\ell}}\n{\phi}_{L_{\ell}^{p'}}\leq
\n{u}_{C^{\ell}(M)}.
$$
Il vient
$$
\int_{M} u\omega\wedge\phi=\int_{M} \omega\wedge u\phi\leq
\n{u}_{C^{\ell}(M)}\n{\omega}_{L_{-\ell}^{p}(M)},
$$
donc
$$
\n{u\omega}_{L_{-\ell}^{p}(M)}\leq
\n{u}_{C^{\ell}(M)}\n{\omega}_{L_{-\ell}^{p}(M)}.\qed
$$

\remarque Plus g{\'e}n{\'e}ralement, si $\phi$ est une forme diff{\'e}rentielle
de classe $C^\ell$ sur $M$, le produit ext{\'e}rieur par $\phi$ est un op{\'e}rateur born{\'e} sur $L_{-\ell}^{p}$.

\begin{prop}
\label{composantes}
Soit $G$ un groupe de Lie muni d'une m{\'e}trique riemannienne invariante {\`a} gauche. Une forme diff{\'e}rentielle sur $G$ est dans $L_{-\ell}^{p}$ si et seulement si ses composantes dans un rep{\`e}re de formes invariantes {\`a} gauche sont des fonctions dans $L_{-\ell}^{p}$.
\end{prop}

\preuve
Par invariance {\`a} gauche, il suffit de le v{\'e}rifier pour la boule
unit{\'e} $B$ de $G$. Soit
$\theta_{1},\cdots,\theta_{n}$ une base de $1$-formes invariantes {\`a}
gauche. Soit
$\{i_{1},\cdots,i_{n}\}$ une permutation de $\{1,\cdots,n\}$, soit $u$ une
fonction de
classe $C^\ell$ {\`a} support compact dans $B$. Il existe une constante
in\-d{\'e}\-pen\-dan\-te de $u$
telle que
$$
\con
\n{u}_{L_{\ell}^{p'}(B)}\leq \n{u\theta_{i_{1}}\wedge\cdots\wedge
\theta_{i_{k}}}_{L_{\ell}^{p'}(B)}\leq\con
\n{u}_{L_{\ell}^{p'}(B)}.
$$
Par dualit{\'e}, on en tire, pour toute fonction $v$ de $L_{-\ell}^{p}(B)$
$$
\con
\n{v}_{L_{-\ell}^{p}(B)}\leq\n{v\theta_{i_{k+1}}\wedge\cdots\wedge
\theta_{i_{n}}}_{L_{-\ell}^{p}(B)}\leq\con
\n{v}_{L_{-\ell}^{p}(B)}.
$$
\qed

\remarque En particulier, dans la boule unit{\'e} de $\R^n$, une forme
diff{\'e}rentielle est
dans $L_{\ell}^{p}$ si et seulement si ses composantes, des fonctions, le sont.

\subsection{R{\'e}gularisation}
Dans cette section, on construit des op{\'e}rateurs de r{\'e}gularisation $S$ et $T$ qui vont constituer une homotopie entre formes $L^p$ et formes $L_{-1}^p$.

\medskip

Soit $M$ une vari{\'e}t{\'e} riemannienne compacte sans bord. Notons
${\cal H}$ l'espace des formes harmonique sur $M$. Le laplacien $\Delta$
admet un inverse born{\'e} sur l'orthogonal $L^2$ de ${\cal H}$. On note $\Gamma=\Delta^{-1}\circ\Pi$ la composition avec le projecteur orthogonal sur ${\cal H}^{\perp}$. Comme $\Delta$ commute avec $d$ et $d^*$,
il en est de m{\^e}me de $\Gamma$.

\begin{lemme}
Soit $M$ une vari{\'e}t{\'e} riemannienne compacte sans bord. Soit $p>1$. Les op{\'e}rateurs $T_{M}=d^{*}\Gamma$ et $S_{M}=1-dT_{M}-T_{M}d$ sont born{\'e}s de $L^{p}_{-1}\Omega^{*}(M)$ dans $L_{}^{p}\Omega^{*}(M)$. De plus, $S_M$ est borné de $L^{p}_{-1}\Omega^{*}(M)$ dans $L_{\ell}^{p}\Omega^{*}(M)$ pour tout $\ell$. 
\end{lemme}

\preuve
Comme $\Gamma$ commute avec $d$, $S_{M}=1-dT_{M}-Td_{M}=1-\Delta \Gamma=1-\Pi$ est le projecteur orthogonal sur ${\cal H}$. Il est donc aussi r{\'e}gularisant qu'on veut.

Montrons que $\Gamma$ est born{\'e} de $L^{p'}$ dans $L_{2}^{p'}$. L'estim{\'e}e $L^{p'}$ pour le laplacien, \cite{ADN}, s'{\'e}nonce comme suit. Pour $\phi\in L_{2}^{p'}$,
$$
\n{\phi}_{L_{2}^{p'}}\leq\con(\n{\phi}_{L^{p'}}+\n{\phi}_{L^{p'}}).
$$
Montrons d'abord que l'image par $\Delta$ de ${\cal H}^{\perp}\cap
L_{2}^{p'}$ est ferm{\'e}e dans $L^{p'}$. Soit $\omega\in L^{p'}$, soit $\phi_{j}\in{\cal H}^{\perp}\cap L_{2}^{p'}$ une suite telle que
$\Delta(\phi_{j})$ converge dans $L^{p'}$ vers $\omega$. Alors $\phi_j$ est
born{\'e}e dans $L^{p'}$. Sinon, $\psi_{j}=\phi_{j}/\n{\phi_{j}}_{L^{p'}}$ est born{\'e} dans $L_{2}^{p'}$ donc, par compacit{\'e} du plongement de Sobolev, on peut supposer que $\psi_{j}$ converge fortement dans $L^{p'}$ vers un $\psi$ de norme $1$, et converge faiblement dans ${\cal
H}^{\perp}\cap L_{2}^{p}$. Or $\Delta\psi_{j}$ tend vers $0$, donc
$\psi\in{\cal H}$. Il vient $\psi=0$, contradiction. On conclut que $\phi_j$ est born{\'e}e dans $L^{p'}$. On peut donc supposer que $\phi_j$ converge faiblement dans ${\cal H}^{\perp}\cap L_{2}^{p'}$ vers $\phi$, et $\omega=\Delta\phi$ est dans l'image du laplacien. On conclut que $\Delta({\cal H}^{\perp}\cap L_{2}^{p'})$ est ferm{\'e} dans $L^{p'}$.

Le laplacien {\'e}tant sym{\'e}trique, son image est exactement
${\cal H}^{\perp}\cap L^{p'}$. L'op{\'e}rateur $\Delta$ est une bijection
continue entre les espaces de Banach ${\cal H}^{\perp}\cap L_{2}^{p'}$ et ${\cal
H}^{\perp}\cap L^{p'}$. C'est donc un isomorphisme. On conclut que $\Gamma$ est born{\'e} de $L^{p'}$ dans $L_{2}^{p'}$.

Le m{\^e}me argument montre que, pour tout $\ell\in\Z$, $\Gamma$ est born{\'e} de $L_{\ell}^{p'}$ dans $L_{\ell+2}^{p'}$. En effet, l'estim{\'e}e elliptique est vraie sur $L_{\ell}^{p'}$. En particulier, l'adjoint $L^2$ de $T_{M}$, l'op{\'e}rateur $T^{*}=d\Gamma$ est born{\'e} de $L_{\ell}^{p'}$ dans $L_{\ell+1}^{p'}$.

Par dualit{\'e}, on conclut que $T_{M}$ est born{\'e} de $L_{-\ell}^{p}$ dans $L_{-\ell+1}^{p}$. Pour $\ell=1$, on trouve que $T_{M}$ est born{\'e} de $\lmpo{*}(M)$ dans $\lpo{*}(M)$.\qed

\begin{cor}
Soit $M$ une vari{\'e}t{\'e} riemannienne {\`a} bord compacte. Soit $p>1$. Il existe
des op{\'e}rateurs $S_{M}$ et $T_{M}$ tels que $1=S_{M}+dT_{M}+T_{M}d$ et
pour toutes fonctions lisses $\chi$ et $\chi'$ {\`a} support compact dans l'int{\'e}rieur de $M$, $\chi'\circ S_{M}\circ\chi$, $d\chi'\wedge T_{M}\circ\chi$ et $T_{M}\circ d\chi\wedge\cdot$ sont born{\'e}s de $L^{p}_{-1}\Omega^{*}(M)$ dans
$L_{}^{p}\Omega^{*}(M)$. De plus, $\chi'\circ S_{M}\circ\chi$ est borné de $L^{p}_{-1}\Omega^{*}(M)$ dans $L_{\ell}^{p}\Omega^{*}(M)$.
\end{cor}
\preuve
On peut toujours compl{\'e}ter $M$ en une vari{\'e}t{\'e} compacte sans bord. On utilise ensuite le fait que $L_{-\ell}^{p}$ est un $C^{\ell'}$-module pour tout $\ell'\geq\ell$ (lemme \ref{module} et la remarque qui le suit).\qed

\begin{prop}
\label{regul}
Soit $M$ une vari{\'e}t{\'e} riemannienne {\`a} g{\'e}om{\'e}trie born{\'e}e. Soit $p>1$. Il existe des op{\'e}rateurs born{\'e}s $S$ et $T$ de $\lmpo{*}(M)$ dans $\lpo{*}(M)$, tels que $1=S+dT+Td$. De plus, pour tout $\ell$, $S$ est borné de $L^{p}_{-1}\Omega^{*}(M)$ dans $L_{\ell}^{p}\Omega^{*}(M)$.
\end{prop}

\preuve
En transportant au moyen de diff{\'e}omorphismes contr{\^o}l{\'e}s en norme $C^1$ les op{\'e}rateurs $S_{B}$ et $T_{B}$ de la boule unit{\'e} $B$ de $\R^n$, on obtient pour chaque $z\in M$ des op{\'e}rateurs $S_z$ et $T_z$ sur la boule $B(z,\epsilon)$ de $M$ tels que $\chi'_{z}\circ S_{z}\circ\chi_{z}$,
$d\chi'_{z}\wedge T_{z}\circ\chi_{z}$ et $T_{z}\circ d\chi_{z}\wedge\cdot$
soient born{\'e}s
de $L^{p}_{-1}\Omega^{*}(B(z,\epsilon))$ dans
$L_{}^{p}\Omega^{*}(B(z,\epsilon))$ en fonction seulement des normes
$C^1$ de
$\chi'_{z}$ et $\chi_{z}$. On choisit les fonctions $\chi'_{z}$ et
$\chi_{z}$ de sorte que
pour tout $x\in M$, $\int \chi'_{z}\chi_{z}(x)\,dz=1$.

On pose, pour $\omega\in L_{-1}^{p}(M)$,
$$
T\omega=\int_{M} \chi'_{z}T_{z}(\chi_{z}\omega_{\mid B(z,\epsilon)})\,dz.
$$
Pour chaque $z\in M$, $\chi_{z}\omega_{\mid B(z,\epsilon)}\in
L_{-1}^{p}(B(z,\epsilon))$ et
$$
\n{\chi'_{z}T_{z}(\chi_{z}\omega_{\mid
B(z,\epsilon)})}_{L_{}^{p}(B(z,\epsilon))}\leq\con
\n{\omega_{\mid B(z,\epsilon)}}_{L_{-1}^{p}(B(z,\epsilon))}
$$
donc, par d{\'e}finition de la norme $L_{-1}^{p}$,
$$
\n{T\omega}_{\lpo{*}}\leq\con\n{\omega}_{\lmpo{*}}.
$$
On calcule
\begin{eqnarray*}
&&d(\chi'_{z}T_{z}(\chi_{z}\omega))+\chi'_{z}T_{z}(\chi_{z}d\omega)\\
&=&d\chi'_{z}\wedge
T_{z}(\chi_{z}\omega)+\chi'_{z}(T_{z}d+dT_{z})(\chi_{z}\omega)-\chi'_{z}T_{z
}(d\chi_{z}\wedge\omega)\\
&=&\chi'_{z}(1-S_{z})(\chi_{z}\omega)+d\chi'_{z}\wedge
T_{z}(\chi_{z}\omega)-\chi'_{z}T_{z}(d\chi_{z}\wedge\omega)\\
&=&\chi'_{z}\chi_{z}\omega-{\tilde S}_{z}(\omega)
\end{eqnarray*}
o{\`u}
$$
{\tilde
S}_{z}(\omega)=\chi'_{z}S_{z}(\chi_{z}\omega)+\chi'_{z}T_{z}(d\chi_{z}\wedge
\omega)
-d\chi'_{z}\wedge T_{z}(\chi_{z}\omega)
$$
donc ${\tilde S}_{z}$ est born{\'e} de
$L_{-1}^{p}(B(z,\epsilon))$
dans $L_{}^{p}(B(z,\epsilon))$, uniform{\'e}ment en $z$.
On pose
$$
S=\int_{M} {\tilde S}_{z}\,dz.
$$
Alors $1=S+dT+Td$, et $S$ et $T$ sont born{\'e}s de $\lmpo{*}(M)$ dans $\lpo{*}(M)$.\qed

\subsection{Comparaison des normes \texorpdfstring{$\psip{*}$}{}
 et \texorpdfstring{$\omp{*}$}{}
}

\begin{defi}
On note $\Omega_{-\ell}^{*,p}(M)$ l'espace des formes
diff{\'e}rentielles qui sont dans $L_{-\ell}^{p}\Omega^{*}(M)$ ainsi que leur diff{\'e}rentielle ext{\'e}rieure, muni de la norme
$$
\n{\omega}_{\Omega_{-\ell}^{*,p}}=(\n{\omega}_{L_{-\ell}^{p}\Omega^{*}(M)}^{p}+
\n{d\omega}^{p}_{L_{-\ell}^{p}\Omega^{*}(M)})^{1/p}.
$$
\end{defi}

\begin{prop}
\label{l-1p=lp+dlp}
Soit $M$ une vari{\'e}t{\'e} riemannienne de dimension $n$ {\`a} g{\'e}om{\'e}trie born{\'e}e. Alors les normes de $\psip{*}(M)$ et de $\omp{*}(M)$ sont {\'e}quivalentes.
\end{prop}
\preuve
Par construction, $\psip{*}\subset \lmpo{*}$ et $\psip{*}$ est
stable par $d$ donc $\psip{*}\subset \Omega_{-1}^{*,p}$.

Inversement, soit $\omega$ une forme dans $L_{-1}^{p}\Omega^{*}$ telle que $d\omega\in L_{-1}^{p}\Omega^{*}$. Utilisons les op{\'e}rateurs de r{\'e}gularisation de la proposition \ref{regul}. Alors
$\omega=(S+Td)\omega+dT\omega$ o{\`u} $(S+Td)\omega\in \lpo{*}$ et
$T\omega\in \lpo{*}$, donc $\omega\in \psip{*}$.\qed

\begin{cor}
\label{comp}
Soit $G$ un groupe de Lie muni d'une m{\'e}trique riemannienne invariante {\`a} gauche. Une forme diff{\'e}rentielle $\omega$ sur $G$ est dans $\psip{*}(G)$ si et seulement si les composantes de $\omega$ et de $d\omega$
dans un rep{\`e}re de formes invariantes {\`a} gauche sont des fonctions dans
$L_{-1}^{p}(G)$.
\end{cor}

\preuve
Combiner les propositions \ref{composantes} et \ref{l-1p=lp+dlp}.\qed

\subsection{Application aux formes basiques sur les produits semi-directs}

\begin{lemme}
\label{compbas}
Soit $G:\R\ltimes_{\alpha} H$ un groupe de Lie produit semi-direct muni d'une m{\'e}trique riemannienne invariante {\`a} gauche. On suppose que $\alpha$ est semi-simple. Soit $W=\mathrm{sp}(\L^{*}\alpha^{\top})$ l'ensemble des valeurs propres de l'action de $\alpha$ sur les formes. Pour chaque $w\in W$, il existe un sous-espace vectoriel invariant à gauche $\mathcal{F}^{w,p}$ de l'espace des distributions sur $H$ ayant la propriété suivante. Soit $(\psi_{j})$ une base de $\L^{*}\H\otimes\mathbb{C}$ formée de vecteurs propres (relativement à des valeurs propres $w_j$) de $\L^{*}\alpha^{\top}$. Soit $\omega$ une forme différentielle sur $H$, soient 
\begin{eqnarray*}
\omega=\sum_{j}u_{j}\psi_{j}, \quad d\omega=\sum_{j}v_{j}\psi_{j}.
\end{eqnarray*}
les composantes de $\omega$ et de sa différentielle dans une base de formes invariantes à gauche. Alors $\omega\in\bp{*}(G,\xi)$ si et seulement si, pour tout $j$, $u_{j}$ et $v_{j}\in\mathcal{F}^{w_{j},p}$.
\end{lemme}

\preuve
Soit $\psi\in\L^{k}\H^{*}$ un vecteur propre complexe de $\L^{k}\alpha^{\top}$, $\L^{k}\alpha^{\top}(\psi)=w\psi$. Considérons $\psi$ comme un élément de $\L^{k}\G^{*}\otimes\mathbb{C}$, et par suite, comme une forme invariante à gauche sur $G$. Alors 
\begin{eqnarray*}
\pi^{*}(\psi_{|H})=e^{-tw}\psi.
\end{eqnarray*}
En effet, le champ $\xi=\ddt$ est invariant à gauche, son flot $\phi_t$ est un groupe à un paramètre de translations à droite, 
\begin{eqnarray*}
\phi_t =R_{\exp(t\xi)}=L_{\exp(t\xi)}\circ Ad_{\exp(t\xi)}
\end{eqnarray*}
La forme $\phi_{t}^{*}\psi$ est invariante à gauche. Sa valeur à l'origine est $(Ad_{\exp(t\xi)})^{*}\phi$. Or
\begin{eqnarray*}
 Ad_{\exp(t\xi)}=\exp(t\mathrm{ad}_{\xi})=\exp(t\alpha),
\end{eqnarray*}
d'où
\begin{eqnarray*}
\phi_{t}^{*}\psi=\exp(t\L^{k}\alpha^{\top})(\psi)=e^{tw}\psi.
\end{eqnarray*}
Au point $(t,h)$ de $G$, $\pi^{*}(\psi_{|H})=\phi_{-t}^{*}\psi_{|H}=e^{-tw}\psi$.

Soit $(\psi_{j})$ une base de $\L^{k}\H^{*}\otimes\mathbb{C}$ formée de vecteurs propres de $\L^{k}\alpha^{\top}$, $\L^{k}\alpha^{\top}(\psi_{j})=w_{j}\psi_{j}$. Alors $\omega=\sum_{j}u_{j}\psi_{j,|H}$ est l'écriture de $\omega$ dans une base de formes invariantes à gauche sur $H$, et
\begin{eqnarray*}
\pi^{*}\omega=\sum_{j}u_{j}e^{-tw_{j}}\psi_{j}
\end{eqnarray*}
est l'écriture de $\pi^{*}\omega$ dans une base de formes invariantes à gauche sur $G$. Par conséquent,
\begin{eqnarray*}
\n{\pi^{*}\omega}_{L_{-1}^{p}(G)}\sim \sum_{j}\n{u_{j}e^{-tw_{j}}}_{L_{-1}^{p}(G)}.
\end{eqnarray*}
On définit $\mathcal{F}^{w,p}$ comme l'espace des distributions $u$ sur $H$ telles que la norme $\n{u_{j}e^{-tw_{j}}}_{L_{-1}^{p}(G)}$ est finie. Alors $\pi^{*}\omega\in L_{-1}^{p}(G)$ si et seulement si pour tout $j$, $u_{j}\in \mathcal{F}^{w_{j},p}$. Par la définition de $\bp{*}(G,\xi)$ vu comme espace de formes différentielles sur $H$ et le Corollaire \ref{comp}, 
\begin{eqnarray*}
\omega\in \bp{*}\quad\Leftrightarrow\quad \forall j,~u_{j}\textrm{ et }v_{j}\in \mathcal{F}^{w_{j},p}.\qed
\end{eqnarray*}

\remarque
Dans \cite{Pweb}, l'espace $\mathcal{F}^{w,p}$ est décrit dans quelques cas particuliers. Par exemple, si $H=\R^{n-1}$ et $\alpha$ est l'identité, il n'y a qu'un poids $w$ dans chaque degré $k$, c'est $w=k$. Dans ce cas, $\mathcal{F}^{w,p}$ est l'espace de Besov $\d B_{p,p}^{-w+(\frac{n-1}{p})}$ lorsque celui-ci a un sens, et est nul sinon. On verra en section \ref{res} que $\mathcal{F}^{w,p}$ est souvent nul.

\begin{cor}
\label{bornef}
Soit $G:\R\ltimes_{\alpha} H$ un groupe de Lie produit semi-direct muni d'une m{\'e}trique riemannienne invariante {\`a} gauche. On suppose que $\alpha$ est semi-simple. Soit $L\in \mathrm{End}(\L^{*}\H^*)$ un endomorphisme qui commute avec $\L^{*}\alpha^{\top}$. Soit $\omega\in\bp{*}(G,\xi)$, vue comme forme différentielle sur $H$. Si $dL\omega\in\bp{*}(G,\xi)$, alors $L\omega\in\bp{*}(G,\xi)$.
\end{cor}

\preuve
On applique le Lemme \ref{compbas}. Soit $\omega=\sum_{w\in W}\omega_{w}$. Par hypothèse, les composantes de $\omega_{w}$ appartiennent à $\mathcal{F}^{w,p}$. Alors $L(\omega)=\sum_{w\in W}L(\omega_{w})$, et les composantes de $L(\omega_{w})$, qui sont des combinaisons linéaires à coefficients constants de celles de $\omega_{w}$, appartiennent elles aussi à $\mathcal{F}^{w,p}$. Par hypothèse, il en est de même pour la forme $dL\omega$. On conclut que $\omega\in\bp{*}(G,\xi)$.\qed

\subsection{Equivalence d'homotopie}

\begin{prop}
\label{regularisation}
Soit $M$ une vari{\'e}t{\'e} riemannienne à géométrie bornée. L'injection des complexes $\op{*}(M)\subset \psip{*}(M)$ est une équivalence d'homotopie.
\end{prop}

\preuve
Soit $I$ l'injection $I:\op{*}(M)\to\omp{*}(M)$. Soient $S$ et $T$ les op{\'e}rateurs fournis par la proposition \ref{regul}. La relation $S=1-dT-Td$ entra{\^\i}ne que $dS=Sd$, donc $S$ envoie $\omp{*}$ dans $\op{*}$. Par conséquent, $SI=S$ et $IS=S$ sont des morphismes de complexes. Comme $S$ et $T$ sont bornés de $\lpo{*}$ dans lui-même, et comme $dT=1-S-Td$, $T$ est borné de $\op{*}$ dans lui-même. Comme $Td=1-S-dT$, $T$ est borné de $\psip{*}$ dans lui-même. Comme, $1-SI=dT+Td$ sur $\op{*}$ et $1-IS=dT+Td$ sur $\omp{*}$, on a bien affaire à une homotopie de complexes.\qed

\section{Propriétés algébriques des formes basiques}
\label{res}

A l'exception du degr{\'e} $1$ (voir \cite{P1}), ou du cas de
l'espace hyperbolique r{\'e}el (voir \cite{Pweb}), on ne sait pas
calculer compl{\`e}tement les espaces $\bp{k}(M)$. N{\'e}anmoins, on va montrer que dans certains cas, certaines composantes d'une forme de $\bp{k}(M)$ sont automatiquement nulles

\subsection{Deux m{\'e}canismes}
\label{dis}

\exemple
\label{exemple}
Soit $G=\R\ltimes_{\alpha} \R^2$ le produit semi-direct d{\'e}fini par la
matrice $\d\alpha=\begin{pmatrix}-1&0\\ 0&2\end{pmatrix}$. Notons
$\pi:G\to H$ la projection. Soit $\beta=a\,dx+b\,dy$ une $1$-forme
diff{\'e}rentielle non nulle sur $H$. Alors $\phi=\pi^{*}\beta$ n'est jamais
dans $L^2$. En effet, notons $\beta_{+}=a\,dx$ et $\beta_{-}=b\,dy$. Il
vient
$$
\n{\phi}_{L^{2}(G)}^{2}\sim\int_{\R}(\n{\beta_{+}}_{L^{2}(\R^{2})}^{2}e^{3t}
+\n{\beta_{-}}_{L^{2}(\R^{2})}^{2}e^{-3t})\,dt
$$
donc
$$
\phi\in L^{2}(\R_{+}\times \R^{2})\Rightarrow \beta_{+}=0,\quad
\phi\in L^{2}(\R_{-}\times \R^{2})\Rightarrow \beta_{-}=0.
$$

La condition $\pi^{*}\beta\in(L^{2}+dL^{2})(\R_{\pm}\times\R^{2})$
impose aussi une restriction, plus faible, sur les composantes de $\beta$. 
On {\'e}crit $\pi^{*}\beta=\omega+d\psi$ o{\`u} $\omega$ et $\psi$
sont dans $L^{2}(\R_{+}\times \R^{2})$. 
On d{\'e}compose encore $\omega=a_{t}+dt\wedge b_{t}$ et $\psi=e_{t}$, o{\`u} les
$a_t$ sont des
$1$-formes et les $b_t$ et $e_t$ des fonctions sur $\R^2$. Alors
$\n{a_{+,t}}+\n{e_{t}}$ est dans
$L^{2}(\R_{+},e^{3t}dt)$. Par cons{\'e}quent il existe une suite $t_j$ tendant
vers $-\infty$ telle que
$a_{-,t_{j}}$ et $e_{t_{j}}$ tendent vers $0$.
D'autre part, pour tout $t$,
$$
\beta=a_{t}+de_{t}=a_{+,t}+a_{-,t}+de_{t}
$$
donc $\beta=\lim a_{-,t_{j}}$ dans $L_{-1}^{2}(\R^{2})$, donc $\beta_{+}=0$. 

Supposons maintenant
que $\pi^{*}\beta\in(L^{2}+dL^{2})(\R_{-}\times\R^{2})$. On trouve que 
$\n{a_{-,t}}$ est dans $L^{2}(\R_{-},e^{-3t}dt)$ donc tend vers $0$ (quitte
{\`a} prendre une sous-suite $t_j$). Il vient $\beta=\lim
a_{+,t_{j}}+de_{t_{j}}$ dans $L_{-1}^{2}(\R^{2})$. On conclut seulement que
$\beta_{-}$ est dans l'adh{\'e}rence $L_{-1}^{2}$ de l'image de
$d_{-}=\frac{\partial}{\partial x}$, ce qui ne dit rien.

Bien que cela ne soit pas directement utile ici, montrons comment un second mecanisme permet de d{\'e}duire de cette annulation partielle l'annulation de la cohomologie $L^2$. On sait qu'une $1$-forme $\beta$ sur $\R^2$ telle que
$\pi^{*}\beta\in {\cal B}^{1,2}(G,\xi)$ est n{\'e}cessairement de la forme $b\,dy$. Si de plus $d\beta=0$, alors la fonction $b$ ne d{\'e}pend pas de $x$, i.e. $\beta$ est invariante par translation le long de l'axe des $x$. Or la norme
$L^{2}+dL^{2}$ de $\pi^{*}\beta$ est une int{\'e}grale sur $\R^2$ d'une quantit{\'e} qui est elle aussi invariante par translation le long de l'axe des $x$. On conclut que cette norme est nulle. Comme ${\cal B}^{1,2}(G,\xi)=0$ et $\xi$ est $(0,2)$-Anosov, le th{\'e}or{\`e}me \ref{thm} donne que $H^{1,2}(G)=0$.

\subsection{Restrictions sur les formes basiques}

\begin{defi}
\label{gamma}
Soit $E$ un fibr{\'e} sur une vari{\'e}t{\'e} $M$. On {\em note} $\Gamma(E)$ l'espace
vectoriel topologique des sections de $E$ {\`a} coefficients
distributions. Lorsque $M=G$ est un groupe de Lie et $E$ un fibr{\'e} invariant {\`a}
gauche sur $G$, on confondra parfois dans la notation le fibr{\'e} $E$ et sa fibre
au-dessus de l'{\'e}l{\'e}ment neutre.
\end{defi}

\begin{prop}
\label{distr}
Soit $M$ une vari{\'e}t{\'e} riemannienne compl{\`e}te. Soit $\xi$ un champ
de vecteurs unitaire $(k-1,p)$-Anosov et $(k,p)$-Anosov sur $M$. Notons $\phi_t$ le flot de $\xi$. On suppose qu'il existe une hypersurface $H\subset M$ orthogonale à $\xi$ telle que 
$$
\R\times H \to M,\quad (t,x)\to\phi_t (x)
$$
soit un difféomorphisme. Pour $j=k-1,\,k$, on note $\d \Lambda^{j}_{-}\oplus\Lambda^{j}_{-}$ la d{\'e}composition $(j,p)$-Anosov de $\L^j T^* H=\ker\iota_{\xi}$. Soit $\omega$ une $k$-forme diff{\'e}rentielle sur $M$ annul{\'e}e par $\iota_{\xi}$ et invariante par le flot de $\xi$, {\`a} coefficients dans $L^{p}_{-1}(M)$. Alors, vue comme forme différentielle sur $H$,
$$
\omega\in\overline{\opp{k}+d\opp{k-1}}\cap\overline{\opm{k}+d\opm{k-1}},
$$
o{\`u} l'adh{\'e}rence est prise au sens des distributions. En particulier, avec les notations de la Définition \ref{basique},
$$
\bp{k}(M,\xi)\subset\overline{\opp{k}+d\opp{k-1}}\cap\overline{\opm{k}+d\opm{k-1}}.
$$
\end{prop}

\preuve
Soit $\omega\in\bp{k}(M,\xi)$. Comme $\omega\in\psip{k}(M)$, il existe des formes $L^p$ $\beta$ et $\gamma$ telles que $\omega=\beta+d\gamma$. Par hypothèse, $M=\R\times H$, et on utilise la décomposition des formes différentielles et de la différentielle exterieure sur un produit. Ecrivons $\beta=a+dt\wedge b$ et $\gamma=e+dt\wedge f$ o{\`u} $a$, $b$, $e$ et $f$ sont annul{\'e}es par $\iota_{\xi}$. Alors $\omega=a+d_{H}e$ o{\`u} $d_H$ est la diff{\'e}rentielle ext{\'e}rieure dans la direction du facteur $H$. Comme $\xi$ est $(k-1,p)$-Anosov et $(k,p)$-Anosov, les $k-1$-formes et $k$-formes se d{\'e}composent en $a=a_{+}+a_{-}$ et $e=e_{+}+e_{-}$, de sorte que, lorsque $t$ tend vers $+\infty$, $\phi_{t}^{*}a_{+}$ et $\phi_{t}^{*}e_{+}$ tendent vers $0$ en norme $L^p$. Par cons{\'e}quent 
$$
\omega=\lim_{t\to+\infty} \phi_{t}^{*}a_{-}+d_{\nu}\phi_{t}^{*}e_{-}.
$$
au sens des distributions. Ceci entra{\^\i}ne la convergence au sens des distributions des restrictions {\`a} presque tout translaté $\phi_t (H)$. 
On obtient des formes différentielles $a_t$ et $e_t$ d{\'e}finies sur $H$, {\`a} valeurs dans $\Lambda_{-}$, et telles que $a_{t}+de_{t}$ converge vers la restriction de $\omega$ à $H$ lorsque $t$ tend vers $+\infty$. De m{\^e}me, on construit des formes $a'_t$ et $e'_t$ sur $H$ telles que $a'_{t}+de'_{t}$ converge vers $\omega_{|H}$ lorsque $t$ tend vers $-\infty$.
\qed

\begin{cor}
\label{cordistr}
Soit $H$ une variété riemannienne complète. Soit $M=\R\times H$, munie d'une métrique riemannienne complète telle que $\xi=\frac{\partial}{\partial t}$ soit unitaire et orthogonal à $H$. Soit $p\geq 1$ et $k\geq 1$. On garde les notations de la Proposition \ref{distr}.
\begin{enumerate}
  \item Si $\xi$ est $(k-1,p)$-Anosov, alors 
\begin{eqnarray*}
  d\bp{k}(M,\xi)\subset \overline{d(\Gamma(\L^{k}_{+}))}\cap \overline{d(\Gamma(\L^{k}_{-}))}.
\end{eqnarray*}  
  \item Si $\xi$ est $(k-1,p)$-contractant (resp. dilatant), $d\bp{k}(M,\xi)=0$.
  \item Si $\xi$ est $(k-2,p)$-contractant (resp. dilatant) et $(k-1,p)$-Anosov, alors $\bp{k}(M,\xi)\subset\Gamma(\L^{k}_{+})$ (resp. $\Gamma(\L^{k}_{-})$).
\end{enumerate}
\end{cor}

\preuve
1. Supposons d'abord que $\xi$ est $(k-2,p)$- et $(k-1,p)$-Anosov. D'apr{\`e}s la proposition \ref{distr},
\begin{eqnarray*}
\bp{k-1}\subset\overline{\Gamma(\L^{k-1}_{+})+d\Gamma(\L^{k-2}_{+})},
\end{eqnarray*}
d'où
\begin{eqnarray*}
d\bp{k-1}(M,\xi)\subset d(\overline{\opp{k-1}}),
\end{eqnarray*}
et de même en remplaçant $+$ par $-$. L'hypothèse $\xi$ est $(k-2,p)$-Anosov n'est pas nécessaire si on ne s'intéresse qu'à l'image de la différentielle, celle-ci faisant disparaître les formes de degré $k-2$ de la preuve de la Proposition \ref{distr}.

2. Si $\xi$ est $(k-1,p)$-contractant, $\L^{k-1}_{+}=0$, donc $d\bp{k-1}(M,\xi)=0$. 

3. Si $\xi$ est de plus $(k-2,p)$-contractant, $\L^{k-2}_{+}=0$, donc
\begin{eqnarray*}
\bp{k-1}\subset\overline{\Gamma(\L^{k-1}_{+})}.
\end{eqnarray*}
Comme le projecteur sur $\L^{k-1}_{-}$ est donné, dans une base de formes invariantes à gauche sur $H$, par une matrice constante, il est continu au sens des distributions. Une limite de sections distributions de $\L^{k-1}_{+}$ est encore une section distribution de $\L^{k-1}_{+}$, donc $\bp{k-1}\subset\Gamma(\L^{k-1}_{+})$.\qed

\subsection{Une description de l'espace \texorpdfstring{$\bp{*}$}{}
}

\begin{cor}
\label{distr'}
Soit $H$ une variété riemannienne complète. Soit $M=\R\times H$, munie d'une métrique riemannienne complète telle que $\xi=\frac{\partial}{\partial t}$ soit unitaire et orthogonal à $H$. Soit $p\geq 1$ et $k\geq 1$. On suppose que $\xi$ est $(k-2,p)$ (condition vide si $k=1$) et $(k-1,p)$-Anosov. On note $\beta_{+}$ la composante d'une forme diff{\'e}rentielle $\beta$ de $H$ sur $\L^{*}_{+}$, et $d_{+}\beta=(d\beta)_{+}$. On note $D$ l'opérateur
\begin{eqnarray*}
D:\Omega^{k-1,p}(H)\to \Omega^{k,p}_{-1}(H),\quad e\mapsto D(e)=d_{-}(e_{+})-d_{+}(e_{-}).
\end{eqnarray*}
Vu comme sous-espace de formes différentielles sur $H$, $\bp{k}(M,\xi)$ contient l'ima\-ge de $D$ et est contenu dans l'adhérence $L^{p}_{-1}$ de l'image de $\delta$.
\end{cor}

\preuve
D'après la Proposition \ref{distr}, tout élément $\omega\in\bp{k}(M,\xi)$ est limite $L^{p}_{-1}$ d'une suite $\beta_{j}^{\pm} +d\gamma_{j}^{\pm}$ où $\beta_{j}^{\pm}$ est une section de $\Lambda^{k}_{\pm}$ et $\gamma_{j}^{\pm}$ une section de $\Lambda^{k-1}_{\pm}$. Par conséquent, la composante $\omega_-$ est la limite des $\dmoins\gamma_{j}^{+}$, et $\omega_+$ est la limite des $\dplus\gamma_{j}^{-}$, donc $\omega$ est limite des $\dmoins e_{j,+}-\dplus e_{j,-}$ où $e_{j}=\gamma_{j}^{+}-\gamma_{j}^{-}$. Par conséquent, $\bp{k}(M,\xi)$ est contenu dans l'adhérence dans l'espace $\omp{k}(H)$ des formes différentielles de la forme $\dmoins e_{+}-\dplus e_{-}$.

Réciproquement, soit $e\in\op{k-1}(H)$ une forme diff{\'e}rentielle. Fixons une fonction lisse $\chi$ sur $\R$ telle que $\chi=0$ au voisinage de $-\infty$ et $\chi=1$ au voisinage de $+\infty$. Notons $\pi:M\to H$ la projection sur le second facteur. Montrons que la forme
$$
\omega=d\chi(t)\wedge\pi^{*}e
$$
est dans $\op{k}(G)$, et 
$$
P\omega=\pi^{*}(\dmoins e_{+}-\dplus e_{-}).
$$
Comme $d\chi$ est {\`a} support compact, $\omega$ et $d\omega=-d\chi\wedge \pi^* de$ sont dans $L^p$. On calcule 
$$
B\omega=\chi(t)\pi^{*}e_{-}-(1-\chi(t))\pi^{*}e_{+},$$$$
Bd\omega=\chi(t)\pi^{*}(-\dmoins e)-(1-\chi(t))\pi^{*}(-\dplus e)
$$
et enfin 
\begin{eqnarray*}
P\omega&=&
\chi\pi^{*}(\dmoins e-d e_{-})-(1-\chi)\pi^{*}(\dplus e-de_{+})\\
&=&\pi^{*}(\dmoins e_{+}-\dplus e_{-}).
\end{eqnarray*}
Ceci prouve que $\dmoins e_{+}-\dplus e_{-}\in\bp{k}(M,\xi)$, lorsque cet espace est vu comme un espace de formes différentielles sur $H$.\qed

\subsection{Exemples}

Conservons les notations du Corollaire \ref{distr'}. Pour que $\bp{k}$ soit non nul, il faut que l'un des deux espaces $\L^{k-1}_{+}$ et $\L^{k}_{+}$ soit non nul. Sinon, l'opérateur $D$ est identiquement nul. De même, il faut que l'un des deux espaces $\L^{k-1}_{-}$ et $\L^{k}_{-}$ soit non nul.

\exemple
Pour l'espace hyperbolique réel, pour tout $p\geq 1$ non critique dans les degrés adéquats, il existe au plus un degré $k$ pour lequel $\bp{k}(M,\xi)\not=0$. En effet, comme le montre le tableau,
\begin{center}
\begin{tabular}{|c|ccccccc|}
\hline
 $p$  &$1$&&$\frac{n-1}{k}$&&$\frac{n-1}{k-1}$&&$+\infty$\\
\hline
$\L^{k-1}_{+(p)}$&&$0$&&$0$&&$\L^{k-1}_{k-1}$&\\\hline
$\L^{k}_{+(p)}$&&$0$&&$\L^{k}_{k}$&&$\L^{k}_{k}$&\\\hline\hline
$\L^{k-1}_{-(p)}$&&$\L^{k-1}_{k-1}$&&$\L^{k-1}_{k-1}$&&$0$&\\\hline
$\L^{k}_{-(p)}$&&$\L^{k}_{k}$&&$0$&&$0$&\\\hline
\end{tabular}
\end{center}
$\L^{k-1}_{+}$ et $\L^{k}_{+}$ (ou bien $\L^{k-1}_{-}$ et $\L^{k}_{-}$) sont simultanément nuls, sauf lorsque $\frac{n-1}{k}<p<\frac{n-1}{k-1}$, i.e. lorsque $k$ est la partie entière de $\frac{n-1}{p}+1$.

\exemple
Pour les autres espaces symétriques de rang un, le complexe $\bp{*}$ est de longueur au plus $\f $, il est non nul seulement dans les $\f $ degrés successifs $\frac{\f m+\f -2}{p}+1,\ldots,\frac{\f m+\f -2}{p}+\f $.

\section{Restrictions sur la torsion}
\label{tor}

\subsection{Cas où le champ est \texorpdfstring{$(k-1,p)$}{}-contractant}
\label{torcontr}

\begin{prop}
\label{torsioncontr}
Soit $M$ une vari{\'e}t{\'e} riemannienne. Soit $\xi$ un champ de vecteurs unitaire, complet, $(k-2,p)$-Anosov et $(k-1,p)$-contractant (resp. $(k-1,p)$-dilatant) sur $M$. Alors $\tp{k}(M)=0$.
\end{prop}

\preuve
D'apr{\`e}s la proposition \ref{distr}, $d\bp{k-1}(M,\xi)=0$, car le
sous-fibr{\'e} $\L^{k-1}_{+}$ (resp. $\L^{k-1}_{-}$) est nul. Par cons{\'e}quent, la cohomologie en degr{\'e} $k$ du complexe $\bp{*}(M,\xi)$ est s{\'e}par{\'e}e. Comme remarqué au Corollaire \ref{cordistr}, l'hypothèse que $p$ est non critique en degr{\'e} $k-2$ n'est pas nécessaire pour cela, mais on en a besoin pour appliquer le Th{\'e}or{\`e}me \ref{thm}, qui donne un isomorphisme $\tp{k}(M)=T^{k}(\bp{*}(M,\xi))=0$.\qed 

\exemple
Avec la Proposition \ref{pince}, on retrouve le résultat principal de $\cite{Ppince}$ : si la courbure sectionnelle est suffisamment pincée, la torsion s'annulle.

\exemple
Pour les espaces symétriques de rang un, on trouve que la torsion s'annulle en degré $k$ lorsque $k=1$ (pour tout $p$) et sinon, si $\frac{\f m+\f -2}{p}<\min\sigma(k-1)$ ou bien $\frac{\f m+\f -2}{p}>\max\sigma(k-1)$, ce qui se traduit par $p<\frac{\f m+\f -2}{k-2+\min\{k,\f \}}$ ou $p>\frac{\f m+\f -2}{k-1+\max\{0,k-1-\f m+\f \}}$. 

\exemple
Soit $G=\R\ltimes_{\alpha}\R^2$ où $\alpha=diag(1,-1)$. Ce groupe, parfois appelé $Sol$, est unimodulaire. Aucun exposant n'est critique en degr{\'e} 1. Comme $\L_{+}^{1}=0$, $d\bp{1}(G,\xi)$ est nul, donc la cohomologie du complexe $\bp{*}$ en degr{\'e} $2$ est s{\'e}par{\'e}e, c'est l'espace de Besov ordinaire $B_{p,p}^{2/p}(\R^{2})$. V. Goldshtein et M. Troyanov ont montr{\'e} que $\hp{2}(G)$ est non nul pour tout $p>1$, \cite{GT}. D'apr{\`e}s \cite{P1}, $\hp{1}(G)=0$. Par dualité de Poincaré, la cohomologie r{\'e}duite est nulle en degré 2, donc $\tp{2}(G)\not=0$. Ce n'est pas contradictoire : comme $p$ est critique en degr{\'e} 0, seul l'énoncé 1 du th{\'e}or{\`e}me \ref{thm} s'applique et donne une surjection $H^{2}(\bp{*})\to \hp{2}(G)$. On voit que la condition supplémentaire ($\xi$ est $(k-2,p)$-Anosov) de la Proposition \ref{torsioncontr} est nécessaire pour avoir l'annulation de la torsion.

\subsection{Rôle de la non commutativité}

Dans ce paragraphe, on s'int{\'e}resse aux produits semi-directs
$G=\R\ltimes_{\alpha}H$. 

On va voir que la non commutativit{\'e} du groupe $H$ force parfois la
torsion $\tp{*}(G)$ {\`a} {\^e}tre non nulle sur un intervalle plus grand que celui donné par la Proposition \ref{torsioncontr}.

On identifiera $\bp{*}(G,\xi=\ddt)$ {\`a} un espace de formes dif\-f{\'e}\-ren\-ti\-el\-les sur $H$, {\em not{\'e}} simplement $\bp{k}$. 

\begin{defi}
\label{defdelta}
Si on identifie l'alg{\`e}bre ex\-t{\'e}\-ri\-eu\-re $\L^{*}\H^{*}$ aux formes dif\-f{\'e}\-ren\-tielles invariantes {\`a} gauche sur $H$, alors la dif\-f{\'e}\-ren\-ti\-el\-le ext{\'e}rieure devient un op{\'e}rateur alg{\'e}brique {\em not{\'e}} $\delta:\L^{k-1}\H^{*}\to \L^{k}\H^{*}$.
\end{defi}

Remarquer que $\delta$ commute avec l'action de la dérivation $\alpha$ sur les formes.

\begin{lemme}
\label{K+}
Soit $G=\R\ltimes_{\alpha}H$ un produit semi-direct. On suppose que $p$ est non critique en degr{\'e} $k-1$. On suppose que $\delta$ est injectif sur $\L^{k-1}_{+}$ et que $\delta(\L^{k-1}_{+})\cap(\H^{*}\otimes\L^{k-1}_{+})=\{0\}$. Alors $d\bp{k-1}\subset\bp{k}$ est fermé. Si, de plus, $p$ est non critique en degr{\'e} $k-2$, alors $\tp{k}(G)=0$.
\end{lemme}

\preuve
D'après le Corollaire \ref{cordistr}, si $e\in\bp{k}$ appartient à l'adhérence, de $d\bp{k-1}$, alors $e$ appartient à l'adhérence, au sens des distributions, de $d(\opp{k-1})$. On peut donc écrire $e=\lim da_j$ o{\`u} $a_{j}\in\opp{k-1}$. Soit $(\psi_i)$ une base de $\L^{k-1}\H^*$ formée de vecteurs propres de $\alpha$, et telle que $(\psi_{1},\ldots,\psi_{\ell})$ soit une base de $\L^{k-1}_{+}$. On voit chaque $\psi_i$ comme une forme différentielle invariante à gauche sur $H$. Ecrivons $a_j =\sum_{i=1}^{\ell}f_{ij}\psi_{i}$. Notons $Q:\L^{k}\H^{*}\to \L^{k}\H^{*}/(\H^{*}\otimes\L^{k-1}_{+})$ la projection. Alors, en chaque point,
\begin{eqnarray*}
Q(da_j )&=&Q(\sum_{i=1}^{\ell}df_{ij}\wedge\psi^{i}+\sum_{i=1}^{\ell}f_{ij}\delta\psi^{i})\\
&=&\sum_{i=1}^{\ell}f_{ij}Q\delta\psi^{i}.
\end{eqnarray*}
Par hypothèse, l'op{\'e}rateur alg{\'e}brique $Q\delta$ admet un inverse {\`a} gauche $(Q\delta)^{-1}:\L^{k}\H^{*}/\H^{*}\otimes \L^{k-1}_{+}\to \L^{k-1}_{+}$. On peut supposer que, comme $\delta$ et $Q$, $(Q\delta)^{-1}$ commute avec l'action de $\alpha$ sur les formes. Alors $a_j =(Q\delta)^{-1}Q(da_j)$ converge vers la section distribution $a=(Q\delta)^{-1}Q(e)$ de $\L^{k-1}_{+}$, qui satisfait $da=e$. Sachant que $da$ appartient déjà à $\bp{k-1}(G,\xi)$, le Corollaire \ref{bornef} entraîne que $a$ appartient à $\bp{k-1}(G,\xi)$, et $e$ à $d\bp{k-1}(G,\xi)$.

Si, de plus, $p$ est non critique en degr{\'e} $k-2$, alors le Théorème \ref{thm} s'applique, et $\tp{k}(G)=\tp{k}(\bp{*})=0$.\qed

\subsection{Preuve du Théorème \ref{torsion}}
\label{preuvetorsion}

\begin{prop}
\label{ross}
Soit $M$ un espace sym{\'e}trique de rang $1$ à courbure non cons\-tan\-te, de dimension $n=\f m$ ($\f=2$, $4$ ou $8$). Soit $k=2,\ldots,\frac{\f m+\f-2}{2}$. Notons $\p1{k-1}<\p2{k-1}<\cdots$ les exposants critiques en degr{\'e} $k-1$. Si $p<\p2{k-1}$, alors $\tp{k}(M)=0$ sauf peut-{\^e}tre pour $p=\p1{k-1}$.
\end{prop}

\preuve
Tout espace sym{\'e}trique de rang $1$ {\`a} courbure non constante est isom{\'e}trique {\`a} un produit semi-direct $G=\R\ltimes_{\alpha}N$ o{\`u} $\alpha$ admet les valeurs propres $2$ (de multiplicit{\'e} $\f-1=1$, $3$ ou $7$) et $1$ (de multiplicit{\'e} $\f m-\f$ paire). On note $h=\mathrm{tr}\,\alpha=\f m+\f -2$. En outre, l'espace propre $\N_2$ est central, la diff{\'e}rentielle $\delta$ est nulle sur $\N^{*}_{1}$ et si $e\in\N^{*}_{2}$ est non nul, $\delta(e)\in\L^{2}\N^{*}_{1}$ est non d{\'e}g{\'e}n{\'e}r{\'e}e. L'alg{\`e}bre ext{\'e}rieure se d{\'e}compose en
$$
\L^{k}\N^{*}=\bigoplus_{j=\max\{0,k-\f m+\f \}}^{\min\{k,\f -1\}}\L^{k}_{k+j}\quad\hbox{o{\`u}}\quad
\L^{k}_{k+j}=\L^{k-j}\N_{1}^{*}\otimes\L^{j}\N_{2}^{*},
$$
et $\alpha$ vaut $k+j$ sur $\L^{k}_{k+j}$. On appelle poids l'indice $k+j$. Comme $\delta$ commute avec $\alpha$, il préserve le poids, $\delta(\L^{k}_{k+j})\subset \L^{k+1}_{k+j}$.

Les exposants critiques en degr{\'e} $k$ sont les nombres de la forme $\frac{h}{k+j}$ o{\`u} $j$ d{\'e}crit les entiers compris entre $\max\{0,k-\f m+\f \}$ et $\min\{k,\f -1\}$. Pour $2\leq k\leq \f m-1$, on v{\'e}rifie que 
\begin{eqnarray*}
\p1{k-1}<\p2{k-1}&=&\frac{h}{\min\{2k-2,k+\f -2\}-1}\\
&\leq&\frac{h}{\min\{2k-4,k+\f -3\}}\\
&=&\p1{k-2}<\p2{k-2}.
\end{eqnarray*}
Par cons{\'e}quent, si $p<\p2{k-1}$, alors $p$ est non critique en degr{\'e} $k-2$. Si de plus $p\not=\p1{k-1}$, alors $p$ est non critique en degr{\'e}s $k-1$ et $k-2$.

Si $p<\p1{k-1}$, alors $\L^{k-1}_{+(p)}=0$, et la Proposition \ref{torsioncontr} donne imm{\'e}diatement que $\tp{k}(M)=0$. Supposons $p$ compris entre $\p1{k-1}$ et $\p2{k-1}$. Alors $\L^{k-1}_{+(p)}=\L^{k-1}_{k-1+r}$, où $r=\min\{k-1,\f -1\}$. Lorsque $k\leq \f $, $r=k-1$ donc $\L^{k-1}_{+}=\L^{k}_{2k-2}=\L^{k-1}(\N_{2})^*$. Lorsque $k>\f $, $r=\f -1$ donc $\L^{k-1}_{+}=\L^{k-1}_{k+\f -2}=\L^{k-\f }\N_{1}^* \otimes\L^{\f -1}(\N_{2})^*$. Dans les deux cas, $\delta(\L^{k-1}_{+})$ a pour poids $k-1+r$, alors que dans $\N^* \otimes \L^{k-1}_{+}$, on rencontre seulement les poids $k+r$ et éventuellement $k+r+1$, donc $\delta(\L^{k-1}_{+})\cap(\N^* \otimes \L^{k-1}_{+})=\{0\}$.

Montrons que si $k\leq\frac{h}{2}$, dans les deux cas, $\delta$ est injectif sur $\L^{k-1}_{+}$. 

Supposons d'abord que $k\leq d$. Soit $v_{1},\ldots,v_{\f -1}$ une base du sous-espace $\N_{2}$. Alors les $2$-formes $\delta v_{1},\ldots,\delta v_{\ell}$ sont lin{\'e}airement ind{\'e}pendantes dans $\L^{2}(\N_{1})^{*}$. Si $I$ est un sous-ensemble de $\{1,\ldots,\f -1\}$, notons 
$$
v_{I}=\bigwedge_{i\in I} v_{i}.
$$
Alors 
$$
\delta v_{I}=\sum_{i\in I} \pm \delta v_{i}\otimes v_{I\setminus\{i\}},
$$
et lorsque $I$ d{\'e}crit les sous-ensembles {\`a} $k-1$ {\'e}l{\'e}ments de
$\{1,\ldots,\f -1\}$ et $i$ d{\'e}crit les {\'e}l{\'e}ments de $I$, les {\'e}l{\'e}ments $\delta v_{i}\otimes v_{I\setminus\{i\}}$ sont lin{\'e}airement ind{\'e}pendants dans $\L^{2}(\N_{1})^{*}\otimes\L^{k-2}(\N_{2})^{*}$. Par cons{\'e}quent $\delta$ est injectif sur $\L^{k}_{2k-2}$.

Supposons maintenant que $k>\f $. Pour $i=1,\dots,\f -1$, notons $L_{i}$ le produit ext{\'e}rieur par $\delta(v_{i})$. Notons $[\f -1]=\{1,\ldots,d-1\}$. Alors sur $\L^{k-1}_{k+\f -2}=\L^{k-\f }\N_{1}^* \otimes\L^{\f -1}(\N_{2})^*$, 
$$
\delta(u\otimes v_{[\f -1]})=\sum_{1\leq i\leq \ell} \pm L_{i}(u)\otimes v_{[\f -1]\setminus\{i\}},
$$
donc $\ker\delta=\bigcap_{i}\ker L_{i}$.

Comme les formes $\delta v_{i}$ sont non d{\'e}g{\'e}n{\'e}r{\'e}es, les applications $L_{i}$ sont injectives sur $\L^{q}(\N_{1})^{*}$ tant que
$q\leq \frac{\f m-\f }{2}-1$. On constate que 
$$
k-\f \leq\frac{\f m-\f }{2}-1\Leftrightarrow k\leq\frac{h}{2}.
$$
Par cons{\'e}quent, si $\f <k\leq \frac{h}{2}$, $\delta$ est injectif sur
$\L^{k-1}_{k+\f -2}$.

Du Lemme \ref{K+}, il r{\'e}sulte que si $k\leq \frac{h}{2}$,
alors $\tp{k}(G)=0$.\qed

\remarque
La borne $\frac{h}{2}$ de la proposition \ref{ross} n'est pas optimale.

\subsection{Caractérisation de \texorpdfstring{$d\bp{k-1}$}{}}

\begin{prop}
\label{db}
Soit $M$ un espace sym{\'e}trique de rang $1$ à courbure non cons\-tan\-te, de dimension $n=\f m$ ($\f=2$, $4$ ou $8$). Notons $\tau$ un vecteur directeur de $\L^{\f-1}(\N_{2})^*$. Notons $\p1{\f-1}=\frac{\f m+\f -2}{2\f -2}<\p2{\f-1}=\frac{\f m+\f -2}{2\f -3}<\cdots$ les exposants critiques en degr{\'e} $\f-1$. Supposons $\p1{\f-1}<p<\p2{\f -1}$. Alors
\begin{enumerate}
  \item $\bp{\f}\subset\Gamma(\L^{\f}_{2\f -1}\oplus\L^{\f}_{2\f -2})$.
  \item Soit $\omega=\omega_{2\f -1}+\omega_{2\f -2}\in\bp{\f}$. Alors 
\begin{eqnarray*}
\omega_{2\f-1}=\alpha\wedge \tau,\quad \omega_{2\f-2}=u\delta(\tau)+\psi,
\end{eqnarray*}
où $\alpha\in\Gamma(\N_{1}^*)$, $u$ est une fonction et $\psi$ une $\f$-forme différentielle à valeurs dans un supplémentaire de $\delta(\tau)$.
  \item $\omega\in d\bp{\f -1}$ si et seulement si $\psi=0$ et $\alpha=d_{\N_1}u$ est la différentielle de $u$ restreinte au champ de plans $\N_1$.
\end{enumerate}
\end{prop}

\preuve
Comme $\p1{\f-1}<p<\p2{\f -1}$, $\L^{\f}_{+(p)}=\L^{\f}_{2\f -2}\oplus\L^{\f}_{2\f -1}$ et $\L^{\f-1}_{+(p)}=\L^{\f-1}_{2\f -2}$. De plus, $\N^* \otimes \L^{\f-1}_{2\f -2}=\L^{\f}_{2\f -1}$ (car $(\N_{2})^* \otimes \L^{\f}_{2\f -2}=0$) et $\delta(\L^{\f-1}_{2\f -2})\subset \L^{\f}_{2\f -2}$, donc $d\L^{\f-1}_{+(p)}\subset \Gamma(\L^{\f}_{2\f -2}\oplus\L^{\f}_{2\f -1})$. D'après la Proposition \ref{distr}, $\bp{\f}\subset\overline{\Gamma(\L^{\f}_{2\f -2}\oplus\L^{\f}_{2\f -1})}=\Gamma(\L^{\f}_{2\f -2}\oplus\L^{\f}_{2\f -1})$.

Soit $\beta\in\bp{\f-1}$, $\beta=a\tau$ où $a$ est une fonction. Alors
\begin{eqnarray*}
d\beta=d_{\N_1}a\wedge\tau+a\delta(\tau),
\end{eqnarray*}
où le premier terme est une section de $\L^{\f}_{2\f -1}$ et le second une section de $\L^{\f}_{2\f -2}$. Pour qu'une $\f$-forme $\omega\in\bp{\f}$ appartienne à l'image de $d$, il est nécessaire que sa composante $\omega_{2\f-2}$ soit proportionnelle à $\delta(\tau)$, $\omega_{2\f-2}=u\delta(\tau)$. Alors la seule solution possible de l'équation $d\beta=\omega$ est $\beta=u\tau$, ce qui impose que $\omega_{2\f-1}=d_{\N_1}u\wedge\tau$.

Inversement, supposons que $\omega\in\bp{\f}$ s'écrive $\omega=\alpha\wedge\tau+u\delta(\tau)+\psi$ avec $\psi=0$ et $d_{\N_1}u=\alpha$. Alors $\beta=u\tau$ satisfait $d\beta=\omega$. $d\beta$ appartient à $\bp{\f}$. De plus, $\beta=L\omega$ où $L=\delta^{-1}Q$ et $Q:\L^{\f}\N^* \to\L^{\f}\N^*$ est un projecteur $\alpha$-invariant sur $\delta(\L^{\f-1}_{2\f -2})$. Le Corollaire \ref{bornef} entraîne que $\beta\in\bp{\f-1}$.\qed

\begin{cor}
\label{cordb}
A l'exception du cas du plan hyperbolique complexe ($\f=2$ et $m=2$), une forme $\omega\in\bp{\f}$ fermée appartient à $d\bp{\f-1}$ si et seulement si sa composante $\omega_{2\f-2}$ est proportionnelle à $\delta(\tau)$. 
\end{cor}
 
\preuve
Soit $\omega\in\bp{\f}$ une forme fermée. Ecrivons $\omega=\alpha\wedge\tau+u\delta(\tau)$. Alors 
\begin{eqnarray*}
0&=&d\omega\\
&=&(d\alpha)\wedge\tau+(d_{\N_1}u-\alpha)\wedge\delta(\tau).
\end{eqnarray*}
Si $\f\not=2$ ou $m>2$, alors $\dim(\N_1)>2$, la multiplication extérieure par $\delta(\tau)$ est une somme directe d'opérateurs de multiplication extérieure par des formes symplectiques. Il est donc injectif, on conclut que $\alpha=d_{\N_1}u$, et, d'après la Proposition \ref{db}, $\omega\in d\bp{\f-1}$.\qed

\section{Cup-produit}
\label{cup}

\subsection{Valeur au bord et produit extérieur}

Sous les hypothèses du Théorème \ref{thm}, l'isomorphisme en cohomologie entre le complexe $\Omega^{*,p}$ et le complexe de Besov $\bp{*}$ n'est pas en général compatible avec le produit extérieur des formes différentielles. Toutefois, c'est le cas en degré $k$ en présence d'un champ de vecteurs $(j,p)$-contractant pour tout $j\leq k-1$.

\begin{lemme}
\label{produitext}
Soit $M$ une vari{\'e}t{\'e} riemannienne compl{\`e}te de dimension $n$. Soient $k$, $\ell$ des entiers tels que $k+\ell\leq n$. Soient $p$, $q$, $r$ des réels positifs tels que $\d \frac{1}{p}=\frac{1}{q}+\frac{1}{r}$. Soit $\xi$ un champ de vecteurs unitaire sur $M$, supposé $(j,q)$-contractant pour tout $j\leq k-1$, $(j,r)$-contractant pour tout $j\leq\ell-1$ et $(j,p)$-contractant pour $j= k+\ell-1$. Soit $V$ un ouvert invariant par le flot de $\xi$. Alors l'application $P:\Omega^{*,p}(V)\to\mathcal{B}^{*,p}(V)$, qui est bien définie et induit un isomorphisme en cohomologie en degrés $k$, $\ell$ et $k+\ell$, satisfait, pour $\alpha\in \Omega^{k,q}(V)$, $\beta\in \Omega^{\ell,r}(V)$ fermées, $P(\alpha\wedge\beta)=P(\alpha)\wedge P(\beta)$.
\end{lemme}

\preuve
Par hypothèse, le flot $\phi_t$ de $\xi$ contracte toutes les $i$-formes $L^q$ et toutes les $j$-formes $L^r$ transversalement à $\xi$. Autrement dit, pour un $\eta>0$.
\begin{eqnarray*}
\parallel \Lambda^{i}(d\phi_{t})\parallel\mathrm{det}(d\phi_{t})^{1/q}\leq C\,e^{-\eta t}\quad\textrm{ et }
\parallel \Lambda^{j}(d\phi_{t})\parallel\mathrm{det}(d\phi_{t})^{1/r}\leq C\,e^{-\eta t}.
\end{eqnarray*} 
Comme $\Lambda^{i+j}(\ker(\iota_{\xi}))$ est engendré par les produits extérieurs de formes transverses de degrés $i$ et $j$ respectivement, on en déduit que
\begin{eqnarray*}
\parallel \Lambda^{i+j}(d\phi_{t})\parallel\mathrm{det}(d\phi_{t})^{1/q+1/r}\leq C^2\,e^{-2\eta t},
\end{eqnarray*} 
i.e. que $\xi$ contracte toutes les $i+j$-formes $L^p$ transversalement à lui-même. Mais cela ne donne pas la contraction sur les $k+\ell-1$-formes, il est donc nécessaire de le supposer.

Dans ces conditions, l'application $P$ est, pour les formes fermées de degrés $k$, $\ell$ et $k+\ell$, une valeur au bord, voir Lemme \ref{voisinage}, elle s'écrit
\begin{eqnarray*}
P\omega=\lim_{t\to +\infty}\phi_{t}^{*}(\omega),
\end{eqnarray*}
donc passe au travers des produits extérieurs. Cela prouve que la forme $P(\alpha)\wedge P(\beta)$ est bien définie et qu'elle est égale à $P(\alpha\wedge\beta)$.\qed

\begin{lemme}
\label{produit}
Soit $M$ une vari{\'e}t{\'e} riemannienne compl{\`e}te de dimension $n$. Soit $q>1$. Soit $\xi$ un champ de vecteurs unitaire sur $M$, supposé $(0,q)$-contractant. Soit $V$ un ouvert invariant par le flot de $\xi$. 
\begin{enumerate}
  \item Soit $u$ une fonction sur $V$ telle que $du\in L^q$. Alors $u$ possède une limite le long de presque toute trajectoire de $\xi$, notée $u_{\infty}$. L'application $P:\Omega^{1,q}(V)\to\mathcal{B}^{1,q}(V)$, qui est bien définie et induit un isomorphisme en cohomologie en degré $1$, satisfait 
$$
P(du)=d(u_{\infty}).
$$
  \item Si $u$ et $v$ sont deux fonctions bornées à différentielles $L^q$ sur $V$, alors $d(uv)\in L^q$ et
\begin{eqnarray*}
P(d(uv))=u_{\infty}P(dv)+v_{\infty}P(du).
\end{eqnarray*}
\end{enumerate}
\end{lemme}

\preuve
1. Comme $du$ est $L^q$ et fermée et $\xi$ est $(0,q)$-contractant, la fonction $u-Bdu$ est bien définie sur $V$. Sa différentielle au sens des distributions vaut 
\begin{eqnarray*}
d(u-Bdu)=(1-dB)(du)=P(du),
\end{eqnarray*}
elle appartient à $\mathcal{B}^{1,q}(V,\xi)$, donc $u-Bdu$ est constante sur les trajectoires de $\xi$. Par définition de l'opérateur $B$,
\begin{eqnarray*}
Bdu&=&-\int_{0}^{+\infty}\phi_{t}^{*}(\iota_{\xi}du)\,dt\\
&=&-\int_{0}^{+\infty}(\xi u)\circ\phi_{t}\,dt\\
&=&u-\lim_{t\to +\infty}u\circ\phi_{t},
\end{eqnarray*}
donc $u-Bdu=\lim_{t\to +\infty}u\circ\phi_{t}$, qu'on peut noter $u_{\infty}$.

2. La différentielle du produit $d(uv)=udv+vdu$ est dans $L^q$. Il vient
\begin{eqnarray*}
P(d(uv))&=&d((uv)_{\infty})\\
&=&u_{\infty}dv_{\infty}+v_{\infty}du_{\infty}\\
&=&u_{\infty}P(dv)+v_{\infty}P(du).\qed
\end{eqnarray*}

\subsection{Localisation en courbure négative}

Dans une variété à courbure négative, on peut enrichir la notion de cohomologie $L^p$ (resp. d'algèbre de Royden) en la localisant près du bord à l'infini.

\begin{defi}
\label{locneg}
Soit $M$ une variété riemannienne complète, simplement connexe, à courbure
sectionnelle négativement pincée. On note $\partial_{\infty}M$ son bord à
l'infini. On note $\bar{M}=M\cup\partial_{\infty}M$ la compactification
visuelle de $M$. Soit $\eta$ un point de $\partial_{\infty}M$, soit $\xi$ le champ
de vecteurs de Busemann correspondant. On note $\Omega^{k,p}_{loc}(M,\eta)$
plutôt que $\Omega^{k,p}_{loc}(M,\xi)$ la cohomologie des formes qui sont
$\Omega^{p}$ dans tout ouvert dont l'adhérence dans $\bar{M}$ ne contient pas $\eta$.
\end{defi}

\begin{defi}
\label{locroy}
Soit $M$ une variété riemannienne complète, simplement connexe, à courbure
sectionnelle négativement pincée. Soit $\eta$ un point de
$\partial_{\infty}M$. On note $\mathcal{R}^{q}_{loc}(M,\eta)$
l'algèbre des fonctions sur $M$ qui sont bornées et à différentielle
$L^{q}$ dans tout ouvert dont l'adhérence dans $\bar{M}$ ne contient pas
$\eta$. On note $\bar{\mathcal{R}}^{q}_{loc}(M,\eta)$ le quotient de
$\mathcal{R}^{q}_{loc}(M,\eta)$ par l'idéal des fonctions qui sont bornées
et $L^{q}$ dans tout ouvert dont l'adhérence dans $\bar{M}$ ne contient pas
$\eta$.
\end{defi}

\subsection{Preuve du Th{\'e}or{\`e}me \ref{subalgebra}}

La proposition suivante est une version localisée, plus forte, du Th{\'e}or{\`e}me \ref{subalgebra}.

\begin{prop}
\label{preuvesubalgebra}
Soit $M$ une vari{\'e}t{\'e} riemannienne compl{\`e}te de dimension $n$,
simplement connexe, dont la courbure sectionnelle $K$ satisfait $-1\leq K\leq\delta<0$. Soit $\eta$ un point de $\partial^{\infty}M$. Soient $r>p\geq1$, soit $q$ le réel tel que $\frac{1}{p}=\frac{1}{q}+\frac{1}{r}$. Supposons que
\begin{eqnarray*}
p<1+\frac{n-k}{k-1}\sqrt{-\delta}\quad\mathrm{et}\quad r<1+\frac{n-k+1}{k-2}\sqrt{-\delta}.
\end{eqnarray*}
\begin{enumerate}
  \item Soit $\kappa\in H^{k-1,r}_{loc}(M,\eta)$. Alors l'ensemble 
\begin{eqnarray*}
\mathcal{R}_{\kappa}=\{u\in\mathcal{R}^q (M,\eta)\,|\,[du]\smile\kappa=0\}
\end{eqnarray*}
est une sous-algèbre de $\mathcal{R}^q (M,\eta)$.
  \item Soient $\kappa_{1},\ldots,\kappa_{\ell}\in H^{k-1,r}_{loc}(M,\eta)$. Notons
\begin{eqnarray*}
\mathcal{R}_{\kappa_{1},\ldots,\kappa_{\ell}}=\bigcap_{j=1}^{\ell}\mathcal{R}_{\kappa_{j}}.
\end{eqnarray*}
C'est une sous-algèbre de $\mathcal{R}^q (M,\eta)$. Alors l'ensemble 
\begin{eqnarray*}
\mathcal{M}_{\kappa_{1},\ldots,\kappa_{\ell}}=\{(u_{1},\ldots,u_{\ell})\in(\mathcal{R}^q (M,\eta))^{\ell}\,|\,\sum_{j=1}^{\ell}[du_{j}]\smile\kappa_{j}=0\}
\end{eqnarray*}
est un $\mathcal{R}^{q}_{\kappa_{1},\ldots,\kappa_{\ell}}$-module.
\end{enumerate}
\end{prop}

\preuve
Soit $\xi$ le champ de vecteur de Busemann dont les trajectoires sont contenues dans $V$. D'apr{\`e}s la proposition \ref{pince}, $\xi$ est $(j,r)$-contractant pour tout degr{\'e} $j\leq k-2$, et $(j,p)$-contractant en degr{\'e} $j= k-1$. Il est automatiquement $(0,q)$ contractant. En particulier, il est $(k-2,p)$-Anosov. Si $u\in\mathcal{R}^{q}(M,\eta)$, $du\in\Omega^{1,q}(M,\eta)$. Soit $\omega$ un représentant de la classe $\kappa$, i.e. une $k-1$-forme fermée sur $M$ qui est $L^r$ sur les ouverts dont l'adhérence dans $\bar{M}$ ne contient pas $\eta$. Du lemme \ref{produitext}, il r{\'e}sulte que $P(du\wedge\omega)=P(du)\wedge P(\omega)$.

1. Soient $u$ et $v\in\mathcal{R}^{q}(V)$ deux fonctions telles que $[du]\smile\kappa=0$ et $[dv]\smile\kappa=0$. D'après le théorème \ref{thm}, $P(du)\wedge P(\omega)=0$ et $P(dv)\wedge P(\omega)=0$. Le lemme \ref{produit} donne 
\begin{eqnarray*}
P(d(uv))\wedge P(\omega)
=u_{\infty}P(dv)\wedge P(\omega)+v_{\infty}P(du)\wedge P(\omega)=0.
\end{eqnarray*}
On applique encore le théorème \ref{thm} pour conclure que $[d(uv)]\smile\kappa=0$. On conclut que $\mathcal{R}^{q}_{\kappa}$ est une sous-algèbre.

2. Soient $\kappa_{1},\ldots,\kappa_{\ell}\in H^{k-1,r}(M,\eta)$. Soit $u\in\mathcal{R}^{q}_{\kappa_{1},\ldots,\kappa_{\ell}}$, i.e. $u\in \mathcal{R}^{q}(M,\eta)$ satisfait
\begin{eqnarray*}
[du]\smile\kappa_{j}=0 \quad \textrm{pour~tout~}j=1,\ldots,\ell.
\end{eqnarray*}
Soit $(u_{1},\ldots,u_{\ell})\in\mathcal{M}_{\kappa_{1},\ldots,\kappa_{\ell}}$, i.e. les $u_j$ sont des éléments de $\mathcal{R}^{q}(M,\eta)$ tels que
\begin{eqnarray*}
\sum_{j=1}^{\ell}[du_{j}]\smile\kappa_{j}=0
\end{eqnarray*}
dans $H^{k,p}(M,\eta)$. Pour chaque $j$, choisissons une forme fermée $L^r$ $\beta_{j}$ représentant $\kappa_j$. Alors, d'après le théorème \ref{thm}, $P(\sum_{j}du_{j}\wedge\beta_{j})=0$. Avec les hypothèses sur $p$ et $r$ et le Lemme \ref{produitext}, cela signifie que
\begin{eqnarray*}
\sum_{j}P(du_{j})\wedge P(\beta_{j})=0.
\end{eqnarray*}
De même, pour tout $j$, $P(du)\wedge P(\beta_{j})=0$ dans $\bp{k}(M,\eta)$. Le lemme \ref{produit} donne
\begin{eqnarray*}
\sum_{j}P(d(uu_{j}))\wedge P(\beta_{j})
&=&u_{\infty}(\sum_{j}P(du_{j})\wedge P(\beta_{j}))+\sum_{j}u_{j,\infty}P(du)\wedge P(\beta_{j})\\
&=&0,
\end{eqnarray*}
d'où $P(\sum_{j}d(uu_{j})\wedge\beta_{j})=0$, et enfin $\sum_{j}[d(uu_{j})]\smile\kappa_{j}=0$. On conclut que $(uu_{1},\ldots,uu_{\ell})\in\mathcal{M}_{\kappa_{1},\ldots,\kappa_{\ell}}$, i.e. que $\mathcal{M}_{\kappa_{1},\ldots,\kappa_{\ell}}$ est un $\mathcal{R}^{q}_{\kappa_{1},\ldots,\kappa_{\ell}}$-module.\qed

\subsection{Preuve du Théorème \ref{nosubalgebra}}

\begin{prop}
\label{nosub}
Soit $M=\mathbb{C}H^2$ le plan hyperbolique complexe. Soit $p$ tel que
$2<p<4$. Soit $q=r=2p$. Soit $\eta$ un point du bord à l'infini de $M$. Il
existe une classe $\kappa\in H^{1,r}_{loc}(M,\eta)$ et une fonction $u\in\mathcal{R}^{q}_{loc}(M,\eta)$ telles que
\begin{enumerate}
  \item $[du]\smile\kappa=0$ dans $\hp{2}_{loc}(M,\eta)$,
  \item $[d(u^2)]\smile\kappa\not=0$ dans $\hp{2}_{loc}(M,\eta)$.
\end{enumerate} 
\end{prop}

\preuve
On identifie le bord à l'infini privé de $\eta$ au groupe d'Heisenberg $N$ de dimension 3. Soient $(x,y,z)$ des coordonnées exponentielles sur $N$, de sorte que $(x,y)$ donnent des coordonnées sur $\R^2 =N/[N,N]$. Alors une base de l'espace des formes différentielles invariantes à gauche est $(dx,dy,\tau)$, où $\N_1 =\mathrm{span}(dx,dy)$ et $\N_2 =\mathrm{span}(\tau)$. On calcule $\tau=dz-xdy$, $\delta(\tau)=-dx\wedge dy$. Le produit semi-direct $G=\R\ltimes_{\alpha}N$ agit simplement transitivement sur $M$, on l'identifie à $M$. On a donc une projection $\pi:M\to N$ dont les fibres sont les géodésiques issues de $\eta$, et une quatrième coordonnée, $t$, sur $M$ (c'est la fonction de Busemann centrée en $\eta$).

On considère les fonctions $u_{\infty}(x,y,z)=x$ et $v_{\infty}(x,y,z)=y$ sur $N$. Soit $\chi$ une fonction lisse sur $\R$, nulle sur $]-\infty,0[$, et valant 1 sur $]1,+\infty[$. On définit deux fonctions sur $M$ par $u=\chi(t)\pi^{*}u_{\infty}$ et $v=\chi(t)\pi^{*}v_{\infty}$. Par construction, $u$ et $v$ s'annulent sur l'horoboule $\{t\leq 0\}$.

Vérifions que $u$ et $v\in\mathcal{R}^{2p}_{loc}(M,\xi)$. Soit $W$ un ouvert de $M$ dont l'adhérence $\bar{W}$ dans $\bar{M}$ ne contient pas $\eta$. Alors $\bar{W}\cap \partial_{\infty}M$ est un compact de $N$. Soit $U$ un ouvert borné à bord lisse de $N$ contenant $\bar{W}$, $V=\R\times U$ le cône de base $U$ et de sommet $\eta$. Alors $W\setminus V$ est relativement compact dans $M$ donc, quitte à agrandir $U$, on peut supposer que $W\subset V$. Il existe alors un $T\in\R$ tel que $W\subset V_T :=]T,+\infty[\times U$.

Comme $2p>4$, $\L^1_{-(2p)}=0$ et $u_{\infty}$ et $du_{\infty}$ appartiennent à $L^{2p}(U)$, donc $du=\pi^{*}u_{\infty}\chi'(t)dt + \chi(t)\pi^{*}du_{\infty}$ est $L^{2p}$ sur $V$. De même pour $v$. Par conséquent, $u_{|V}$, $v_{|V}\in\mathcal{R}^{2p}(V)$. A fortiori, $u_{|W}$, $v_{|W}\in\mathcal{R}^{2p}(W)$. On conclut que $u$ et $v\in\mathcal{R}^{2p}_{loc}(M,\xi)$. En particulier, $du$, $dv\in\Omega^{1,2p}_{loc}(M,\eta)$, et $du\wedge dv\in\op{2}_{loc}(M,\eta)$.

Pour $2<p<4$, $\L^1_{+(p)}=\N_{2}^*$, $\L^1_{-(p)}=\N_{1}^*$. Par conséquent, le champ de vecteurs de Busemann $\xi=\ddt$ n'est pas $(1,p)$-contractant. Néanmoins, la conclusion du Lemme \ref{produitext} reste vraie, 
\begin{eqnarray*}
P(du\wedge dv)=P(dv)\wedge P(du).
\end{eqnarray*}
En effet, comme $du_{\infty}$ et $dv_{\infty}$ sont des sections de $\L^1_{-(p)}$, l'opérateur $P$ se comporte comme une valeur au bord. Voici les détails.
\begin{eqnarray*}
du\wedge dv=\chi'(t)\chi(t)dt\wedge\pi^{*}(u_{\infty}dv_{\infty}-v_{\infty}du_{\infty})+\chi(t)^2 \pi^* (du_{\infty}\wedge dv_{\infty}),
\end{eqnarray*}
\begin{eqnarray*}
\iota_{\xi}(du\wedge dv)=\chi'(t)\chi(t)\pi^{*}(u_{\infty}dv_{\infty}-v_{\infty}du_{\infty})
\end{eqnarray*}
est une section de $\L^1_{-(p)}$,
\begin{eqnarray*}
B(du\wedge dv)&=&-(\int_{0}^{+\infty}\chi'(t+s)\chi(t+s)\,ds)\pi^{*}(u_{\infty}dv_{\infty}-v_{\infty}du_{\infty})\\
&=&\frac{1}{2}(\chi(t)^2 -1)\pi^{*}(u_{\infty}dv_{\infty}-v_{\infty}du_{\infty}),
\end{eqnarray*}
dont la différentielle vaut
\begin{eqnarray*}
\chi'(t)\chi(t)dt\wedge\pi^{*}(u_{\infty}dv_{\infty}-v_{\infty}du_{\infty})+(\chi(t)^2 -1)\pi^{*}(du_{\infty}\wedge dv_{\infty}).
\end{eqnarray*}
Il vient
\begin{eqnarray*}
P(du\wedge dv)&=&\chi(t)^2 \pi^* (du_{\infty}\wedge dv_{\infty})+(1-\chi(t)^2)\pi^{*}(du_{\infty}\wedge dv_{\infty})\\
&=&\pi^{*}(du_{\infty}\wedge dv_{\infty})\\
&=&-\pi^{*}d\tau\\
&=&-d\pi^{*}(\tau).
\end{eqnarray*}

Comme $\L^0_{+(p)}=0$, $\mathcal{B}^{1,p}(V,\xi)\subset\Gamma(\N_2^*)$. En fait, la fonction $z$ est une section $L^p$, à différentielle $L^p$, de $\L^0_{-(p)}$ sur $U$. Alors $\dplus z=dz-xdy=\tau$. D'après le Corollaire \ref{distr'}, la forme $\pi^{*}\tau$ appartient à $\bp{1}(V,\xi)$. Cela entraîne que $[P(du\wedge dv)]=0$ dans $H^{2}(\bp{*}(V))$. Le Corollaire \ref{corthm} donne donc $[du\wedge dv]=0$ dans $\hp{2}(W)$. On conclut que $[du]\smile[dv]=0$ dans $\hp{2}_{loc}(M,\eta)$.

Pour la même raison, 
\begin{eqnarray*}
P(d(u^2)\wedge dv)=P(d(u^2))\wedge P(dv)=\pi^{*}(d(u_{\infty}^2)\wedge dv_{\infty}).
\end{eqnarray*}
D'après la Proposition \ref{db}, la forme $\pi^{*}(d(v_{\infty}^2)\wedge dv_{\infty})=\pi^{*}(2xdx\wedge dy)$ n'appartient pas à $d\mathcal{B}^{1,p}(V)$. En effet, elle s'écrit $-2x\delta(\tau)$ mais ne possède pas de composante $d(-2x)\wedge\tau$. Le Lemme \ref{ouvert} donne donc $[d(u^2)\wedge dv]\not=0$ dans $\hp{2}(V_T)$. Comme $W\subset V_T$, $[d(u^2)]\smile[dv]\not=0$ dans $\hp{2}(W)$. On conclut que $[du]\smile[dv]\not=0$ dans $\hp{2}_{loc}(M,\eta)$.\qed

\section{Intervention des quasiisométries}
\label{interqi}

\subsection{Invariance de la cohomologie et du cup-produit}

Voici la version technique du Théorème \ref{qi} dont nous avons besoin. 

\begin{prop}
{\em \cite{qi}}
\label{qitech}
Soient $M$, $M'$ des variétés riemanniennes complètes, simplement connexes,
à courbure pincée négativement. Soit $f:M\to M'$ une quasiisométrie. Soit
$\eta\in\partial_{\infty}M$. Alors $f$ induit, pour tout $p>1$, un isomorphisme
$f^{*}:\hp{*}_{loc}(M',f(\eta))\to\hp{*}_{loc}(M,\eta)$ qui préserve
le cup-produit. De plus, $f$ induit, pour tout $q>1$, un isomorphisme d'algèbres $f^{*}:\bar{\mathcal{R}}^{q}_{loc}(M',f(\eta))\to\bar{\mathcal{R}}^{q}_{loc}(M,\eta)$.
\end{prop}

\subsection{Preuve du Corollaire \ref{thmpince}}

Soit $M$ une vari{\'e}t{\'e} riemannienne compl{\`e}te de dimension $n$,
simplement connexe, dont la courbure sectionnelle $K$ satisfait $-1\leq
K\leq\delta<0$. Notons $M'=\mathbb{C}H^2$. Supposons qu'il existe une quasiisométrie $f:M\to \mathbb{C}H^2$.

Soit $p_0 = 1+2\sqrt{-\delta}$. Si $\delta<-\frac{1}{4}$,
$p_0 >2$. Pour $2<p<p_0$, posons $q=r=2p$. Choisissons un point
$\eta\in\partial_{\infty}M$. Le Théorème \ref{subalgebra} (plus précisément, la Proposition \ref{preuvesubalgebra}) entraîne que, dans $M$, si $\kappa\in H^{1,r}_{loc}(M,\eta)$ et $u\in\mathcal{R}_{\kappa}$, alors $u^{2}\in\mathcal{R}_{\kappa}$. Or le Théorème \ref{nosubalgebra} (plus précisément, la Proposition \ref{nosub}) fournit une classe $\kappa'\in H^{1,r}_{loc}(M,\eta')$ et une fonction $u'\in\mathcal{R}_{\kappa}$ telles que $u'^{2}\notin\mathcal{R}_{\kappa'}$, cela contredit le fait que les quasiisométries préservent la cohomologie, y compris dans sa version localisée, ainsi que le cup-produit. On conclut que $M$ n'est pas quasiisométrique au groupe de Heisenberg.\qed

\vskip1cm

Mots cl{\'e} : Cohomologie $L^{p}$, courbure n{\'e}gative, espace symétrique, espace de Besov, algèbre de Royden. 

Keywords : $L^{p}$-cohomology, negative curvature, symmetric space, Besov space, Royden algebra.

Mathematics Subject Classification :
43A15, 
43A80, 
46E35, 
53C20, 
53C30, 
58A14. 

\vskip 1cm

\noindent Laboratoire de Math{\'e}matique d'Orsay\\
UMR 8628 du C.N.R.S.\\
Universit{\'e} Paris-Sud XI\\
B{\^a}timent 425\\
91405 Orsay\\
France\\
\smallskip
{\tt\small Pierre.Pansu@math.u-psud.fr\\
http://www.math.u-psud.fr/\~{}pansu}

\end{document}